\documentclass[reqno,12pt]{amsart}

\usepackage{amsmath, amscd, amssymb}
\usepackage{latexsym, amssymb, amsmath, amscd, amsfonts}
\usepackage{color}
\usepackage{tabmac}
\usepackage{xcolor}
\usepackage{fancybox}

\usepackage{tikz}
\newcommand*\circled[1]{\tikz[baseline=(char.base)]{
    \node[shape=circle, draw, inner sep=1pt, 
        minimum height=12pt] (char) {\vphantom{1g}#1};}}

\usepackage{pstricks,pst-node}
\usepackage[all]{xy}

\numberwithin{equation}{section}
\newtheorem{lemma}[equation]{Lemma} 
 \newtheorem{proposition}[equation]{Proposition}
\newtheorem{theorem}[equation]{Theorem}
\newtheorem{corollary}[equation]{Corollary}
\newtheorem*{m-theorem1}{Main Theorem 1}  
\newtheorem*{m-theorem2}{Main Theorem 2}
\newtheorem*{m-lemma}{Main Lemma}  
\newtheorem{definition}[equation]{Definition}
\newtheorem{example}[equation]{Example}
\newtheorem{notation}[equation]{Notation}
\newtheorem{remark}[equation]{Remark}

\newcommand{\C}{{\mathbb C}}

\newcommand{\Z}{{\mathbb Z}}

\newcommand{\gl}{{\mathfrak{gl}}}
\newcommand{\GL}{{ \mathrm{GL}} }
\newcommand{\Gln}{{\GL_n}}
\newcommand{\Glp}{\GL_p}

\newcommand{\rM}{\mathrm{M}}

\newcommand{\calr}{\mathcal{R}}
\newcommand{\calq}{\mathcal{Q}}

\newcommand{\Hom}{\mathrm{Hom}}
\newcommand{\LM}{\mathrm{LM}}
\newcommand{\bB}{\mathbf{B}}
\newcommand{\SST}{\mathrm{SST}}
\newcommand{\sgn}{\mathrm{sgn}}

\newcommand{\ve}{\varepsilon}
\newcommand{\calb}{\mathcal{B}}
\newcommand{\calm}{\mathcal{M}}

  \newcommand{\calt}{\mathcal{T}}

\newcommand{\bigs}{\mathbb{S}}
\newcommand{\fgl}{\mathfrak{gl}}
\newcommand{\fm}{\mathfrak{m}}
\newcommand{\fkh}{\mathfrak{h}}

  \newcommand{\fh}{\mathfrak{h}}
  \newcommand{\fu}{\mathfrak{u}}
   
   \newcommand{\fg}{\mathfrak{g}}
   \newcommand{\fb}{\mathfrak{b}}
 \newcommand{\fl}{\mathfrak{l}}

\newcommand{\ot}{\otimes}
\newcommand{\lra}{\longrightarrow}
\newcommand{\ovl}{\overline}

\newcommand{\id}{{\rm id}}
\newcommand{\Sym}{{\mathrm{Sym}}}
\newcommand{\End}{{\mathrm{End}}}
\newcommand{\BS}{\mathbb{S}}

\newcommand{\Lam}{\Lambda}

\newcommand{\beq}{\begin{equation}}
\newcommand{\eeq}{\end{equation}}

\newcommand{\bmtx}{\begin{pmatrix}}
\newcommand{\emtx}{\end{pmatrix}}

\newcommand{\baln}{\begin{aligned}}
\newcommand{\ealn}{\end{aligned}}

\newcommand{\CR}{\mathcal{R}}

 \newcommand{\fgpp}{\fg^{\prime\prime}}
 
 \newcommand{\fq}{\mathfrak{q}}

\begin{document}

\title[Branching algebras]{Branching algebras for  the general linear \\ Lie superalgebra}

\author{Soo Teck Lee}
\address{Department of Mathematics, National University of Singapore,
Block S17, 10 Lower Kent Ridge Road, Singapore 119076, Republic of Singapore}
\email{matleest@nus.edu.sg}

\author {Ruibin Zhang}
\address{School of Mathematics and Statistics,  University of Sydney, Sydney, N.S.W. 2006,
Australia} 
\email{ruibin.zhang@sydney.edu.au}

 
\begin{abstract} 
We develop an algebraic approach to the branching of representations of 
the general linear Lie superalgebra $\fgl_{p|q}(\C)$, 
 by constructing certain super commutative algebras   
whose structure encodes the branching rules.
Using this approach, we derive the branching rules for restricting 
any irreducible polynomial representation $V$  
of $\fgl_{p|q}(\C)$ to a regular subalgebra isomorphic to 
 $\gl_{r|s}(\C)\oplus \gl_{r'|s'}(\C)$, 
$\fgl_{r|s}(\C)\oplus\gl_1(\C)^{r'+s'}$
or $\fgl_{r|s}(\C)$, 
with $r+r'=p$ and $s+s'=q$. 
In the case of $\fgl_{r|s}(\C)\oplus\gl_1(\C)^{r'+s'}$ with $s=0$ or $s=1$ but general $r$, we also construct a basis for the space of $\gl_{r|s}(\C)$ highest weight vectors in $V$; when $r=s=0$, the branching rule leads to explicit expressions for the weight multiplicities of $V$ in terms of Kostka numbers.   
\end{abstract}

\subjclass[2010]{05E10, 15A75, 20G05, 22E46}
\keywords{General linear Lie superalgebras, Reciprocity laws,   Branching rules}

\maketitle


\section{Introduction}

Branching rules and tensor product decompositions are two aspects of the representation theory 
of Lie superalgebras \cite{K}, which are most frequently used in physics, 
particularly in building supersymmetric models of elementary particle (see e.g., \cite{DJ}). 
Much research has been devoted to them 
since the discovery of Lie superalgebras \cite{SNR, CNS} in the 70s, 
but most studies are concerned with special cases with immediate physical applications, 
and general results are rarely obtained.
The generic non-semi-simplicity of finite dimensional representations of Lie superalgebras 
makes the branching rules and tensor product decompositions much more difficult to study 
for Lie superalgebras than for semi-simple Lie algebras \cite{GW}.  
However, polynomial representations of the general linear Lie superalgebra turn out to be an exception. 
They are unitarizable with respect to the compact real form \cite{GZ}, thus are semi-simple. 

We shall  study branching rules for the general linear Lie superalgebra 
by developing an approach inspired by the work of Howe et al.  \cite{HTW1,HTW2}.  
The key feature of this approach is that the branching rules are encoded
in the structure of certain super commutative algebras
to be constructed in this paper. 
 
Let $\fg=\gl_{p|q}$ be a general linear Lie superalgebra (defined precisely in \S \ref{sect:glpq}) and let $\fl$ be a regular Lie super subalgebra  of $\fg$. If $V$ is an  irreducible polynomial representation of $\fg$, then we shall consider the action by $\fl$ on $V$ obtained by restriction. We call a description of how $V$ decomposes into irreducible $\fl$ representations a {\bf\em $(\fg,\fl)$ branching rule.}
We shall determine  $(\fg,\fl)$  branching rules  for the cases when (i) $\fl\cong\gl_{r|s}\oplus \gl_1^{r'+s'}$, (ii) $\fl\cong\gl_{r|s}$, and (iii) $\fl\cong\gl_{r|s}\oplus\gl_{r'|s'}$, where $r,s,r',s'$ are non-negative integers such that $r+r'=p$ and $s+s'=q$, and $\gl_1$ is the one-dimensional Lie algebra over $\C$. In the case when $r=s=0$, the branching rule in (i) leads to explicit expressions for the weight multiplicities in $V$ in terms of Kostka numbers. 
In addition, we also construct a set of linearly independent $\fg$ weight vectors of $V$. It
 forms a basis for the space of $\gl_{r|s}$ highest weight vectors in $V$ if $s=0$ or $s=1$, and 
 coincides with the basis of $V$ obtained in \cite[Theorem 3.3]{CZ} if $r= s=0$.

Let us now briefly describe our approach to solutions of the branching problems.

Fix a positive integer $n$  and consider the natural action by $\gl_n\oplus\gl_{p|q}$ on 
the  $\Z_2$-graded symmetric superalgebra
$\calr:=\bigs(\C^n\otimes\C^{p|q})$ over $\C^n\ot\C^{p|q}$, where $\gl_n=\gl_{n|0}$ is the general linear Lie algebra.
By  the $(\gl_n, \gl_{p|q})$ Howe duality \cite{H1} (also see \cite{CW, Se}), $\calr$ is a multiplicity free sum of irreducible polynomial representations of $\gl_n\oplus\gl_{p|q}$. Let $\fl$ be a subalgebra of $\fg$ of the type (i)-(iii) described above, and let $\fu_n$ and $\fu_\fl$ be the standard maximal nilpotent subalgebras of $\gl_n$ and $\fl$ respectively. We consider the subalgebra 
\[\calq(\fg,\fl)=\calr^{\fu_n\oplus\fu_\fl}\]
 of $\calr$ consisting of all vectors annihilated by the operators from 
$\fu_n\oplus\fu_\fl$. Then $\calq(\fg,\fl)$ is a representation for $\fh_n\oplus \fh_\fg$, where $\fh_n$ and $\fh_\fg$ are Cartan subalgebras of $\gl_n$ and $\fg$ respectively. Thus $\calq(\fg,\fl)$ is a direct sum of eigenspaces for   $\fh_n\oplus \fh_\fg$
\[\calq(\fg,\fl)=\bigoplus_i\calq(\fg,\fl)_i.\]
Each eigenspace $\calq(\fg,\fl)_i$  consists of $\gl_n\oplus\fl$ highest weight vectors of a specific weight, so that its dimension is equal to the multiplicity of an irreducible $\fl$ representation in an irreducible $\fg$ representation. In this sense the structure of $\calq(\fg,\fl)$ encodes a branching rule for $(\fg,\fl)$. In view of this property of $\calq(\fg,\fl)$,  we call it  a {\bf\em branching algebra} for $(\fg,\fl)$.

It turns out that the structure of $\calq(\fg,\fl)$ also encodes a second branching rule for a pair of Lie algebras of the form $(\gl_n^k,\Delta(\gl_n))$,  where $\gl_n^k=\underbrace{\gl_n\oplus \dots\oplus \gl_n}_k$ for any given positive integer $k$ and $\Delta(\gl_n)$ is a copy of $\gl_n$ which embeds in $\gl_n^k$ diagonally. 
The second branching rule is well understood as it is related to either  the Pieri rule  or the Littlewood-Richardson rule which are well known. Hence, we use the second branching rule to deduce the branching rule for $(\fg,\fl)$.

We believe that the approach described above to branching rules of the general linear Lie superalgebra is of general interest and should extend to the other classical series of Lie superalgebras. 

We also construct a linearly independent set of vectors in $\calr$. These vectors are $\gl_n$ highest weight vectors  and are $\fg$ weight vectors. 
 In the case when $s=0$ or $s=1$, these vectors are also highest weight vectors for $\gl_{r|s}$. In fact, we obtain from these vectors a basis for the space of $\gl_{r|s}$ highest weight vectors in an irreducible polynomial $\fg$ representation  $V$  which occurs in $\calr$. We remark that 
 the explicit form of the $\gl_{r|s}$  highest weight vectors  provide 
 more refined information than branching multiplicity.  
 In the case when $r=s=0$, we reproduce the basis for $V$ which was constructed in \cite[Theorem 3.3]{CZ}
 by direct calculation. 
 
Now we relate the present paper to previous work. 
 Branching rules for $\gl_{p|q}$ were extensively investigated in \cite{BR}
by using symmetric group techniques \cite{DJ} through Schur-Weyl duality. 
In particular, a $(\gl_{p|q}, \gl_{r|s})$ branching rule 
for any irreducible polynomial $\gl_{p|q}$-representation  
was given in \cite[Theorem 5.14]{BR}.
A $(\gl_{p|q}, \gl_{r|s}\oplus\gl_{r'|s'}$) branching rule was also obtained in \cite{CLZ} for the oscillator representations by exploring a general principle that
relates branching rules of $\gl_{p|q}$-representations and tensor product decompositions of 
$\gl_n$-representations in the context of $(\gl_n, \gl_{p|q})$ duality. 
We also mention that the same method was applied in \cite{CZ} to obtain a branching rule for the oscillator representations of the orthosymplectic Lie superalgebra, where a $(\gl_{p|q}, \gl_{r|s}\oplus\gl_{r'|s'}$) branching rule for the polynomial representations of $\gl_{p|q}$ was also alluded to (see \cite[Remark 9.4]{CZ}). 
 
The content of this paper is arranged as follows.  In Section 2, we set up notation and also review some basic facts in representation theory which we need in the rest of the paper.  Sections 3 and 4 are devoted to constructing branching algebras for $(\gl_{p|q},\gl_{r|s}\oplus\gl^{r'+s'}_1)$ and for $(\gl_{p|q},\gl_{r|s}\oplus \gl_{r'|s'})$ respectively, and to deriving the associated branching rules. Finally in Section 5, we construct a linearly independent set of vectors in $\calr$ with properties described in the preceding paragraph. 

\medskip
We work over the field $\C$ of complex numbers throughout the paper.
 
\subsection*{Acknowledgement}
The first named author expresses his sincere gratitude to the University of Sydney for
warm hospitality during his visit in   December 2023 - February 2024.

\section{Preliminaries} 
In this section, we will review the {\em iterated Pieri rules} for the general linear Lie algebra $\gl_n$ \cite{HKL, KLW}, and  the Howe duality between $\gl_n$  and the general linear Lie superalgebra $\gl_{p|q}$ \cite{H1} (also see \cite{CLZ, CW, CZ, Se}). These results provide the main technical tools which we will apply to  derive branching rule for polynomial representations of $\gl_{p|q}$. 
Some necessary background material on linear algebra of $\Z_2$-graded vector spaces will also be given.

\subsection{The $(\gl_n, \gl_{p|q})$ Howe duality}
\subsubsection{Linear algebra of $\Z_2$-graded vector spaces}\label{sect:super-facts}
A $\Z_2$-graded vector space $M=M_{\bar 0}\oplus M_{\bar 1}$ is the direct sum of 
the even subspace $M_{\bar 0}$ and the odd subspace $M_{\bar 1}$. 
It is said to be purely even (resp. odd) if $M_{\bar 1}=0$ (resp. $M_{\bar 0}=0$).  
We denote by $[v]=i$ the degree of a homogeneous element $v\in M_{\bar i}$.   
Let $M'=M'_{\bar 0}\oplus M'_{\bar 1}$ be another $\Z_2$-graded vector space. Then the space of linear maps from $M$ to $M'$ is $\Z_2$-graded with 
$\Hom_\C(M, M')= \Hom_\C(M, M')_{\bar 0} \oplus \Hom_\C(M, M')_{\bar 1}$, where any $\varphi\in \Hom_\C(M, M')_{\bar i}$ satisfies $\varphi(M_{\bar j})\subset M'_{\ovl{i+j}}$ for all $i, j=0, 1$. Note in particular that $\End_\C(M)$ has the structure of  a $\Z_2$-graded associative algebra, i.e., an associative superalgebra. 

The usual tensor product of vectros spaces extends to the category of $\Z_2$-graded vector spaces. 
For any objects $M, M'$, we have 
\[
M\ot M'= (M\ot M')_{\bar 0}\oplus (M\ot M')_{\bar 1}, \quad (M\ot M')_{\bar k}=\bigoplus_{\ovl{i+j}=\bar k} M_{\bar i}\ot M' _{\bar j}.
\]
The category of $\Z_2$-graded vector spaces has a canonical symmetry given by  
\beq
\tau_{M, M'}: M\ot M' \lra M'\ot M
\eeq
for any objects $M, M'$, which is defined by the unique linear extension of the map $v\ot v'\mapsto (-1)^{[v][v']} v'\ot v$ for all homogeneous elements $v\in M$ and $v'\in M'$. 
This in particular gives rise to an action of the symmetric group $\Sym_r$ of degree $r$ on the $r$-th tensor power $M^{\ot r}=\underbrace{M\ot M\ot\dots\ot M}_r$ of $M$ for any $r\ge 2$ (see, e.g., \cite{DJ} and \cite{BR}). To describe this, we recall the standard presentation of $\Sym_r$ with generators $s_1, \dots, s_{r-1}$ and defining relations 
\[
s_i^2=s_i,  \quad s_i s_{i+1} s_i = s_{i+1} s_i s_{i+1}, \quad \text{for all valid $i$}. 
\] 
Then $s_i$ acts on $M^{\ot r}$ by 
$
\nu_r(s_i):=\id_M^{\ot (i-1)}\ot \tau_{M, M} \ot \id_M^{\ot (r-i-1)}
$
for all $i$. 

The tensor algebra $T(M)=\oplus_{r\ge 0}M^{\ot r}$ over a $\Z_2$-graded vector space $M$ is a superalgebra, 
which is also $\Z_+$-graded with $M^{\ot r}$ being the degree $r$ homogeneous subspace.
Here, $\Z_+$ denotes the set of all nonnegative integers. Let $J(M)$ be the two-sided ideal of $T(M)$ generated by the elements $v\ot v'- (-1)^{[v][v']} v'\ot v$ for all homogeneous $v, v'\in M$. This is a homogeneous ideal with respect to both the $\Z_2$-grading and $\Z_+$-grading, and hence  
\[
\BS(M):=T(M)/J(M)
\]
is a $\Z_+$-graded associative superalgebra, which is usually referred to as the $\Z_2$-graded symmetric algebra, or supersymmetric algebra, over $M$. It is $\Z_2$-graded commutative in the sense that $x y - (-1)^{[x][y]} y x=0$ for all homogeneous $x, y\in \BS(M)$. 

Note the following obvious facts.
\begin{enumerate}
\item[(i)]  
$
\BS(M)=S(M_{\bar 0})\ot \Lam(M_{\bar 1})
$, where $S(M_{\bar 0})$ is the usual symmetric algebra of $M_{\bar 0}$ and $\Lam(M_{\bar 1})$ is the exterior
algebra of $M_{\bar 1}$; 

\item[(ii)]   $J(M)_r=J(M)\cap M^{\ot r}$ is a $\Sym_r$-submodule, 
and
$M^{\ot r} \cong \BS^r(M) \oplus J(M)_r$ as $\Sym_r$-module  for any  $r\in\Z_+$,  
where $\BS^r(M)$ is the degree $r$ homogeneous subspace of $\BS(M)$;

\item[(iii)]   
$
\BS(M\oplus M')=\BS(M)\ot \BS(M')
$
for any $\Z_2$-graded vector spaces $M, M'$. 

\end{enumerate}
Part (ii) makes use of the semi-simplicity of $M^{\ot r}$ as $\Sym_r$-module.

\subsubsection{The general linear Lie superalgebra $\gl_{p|q}$}\label{sect:glpq}

For any $\Z_2$-graded vector space $V=V_{\bar 0}\oplus V_{\bar 1}$, 
the endomorphism superalgebra $\End_\C(V)$ can be endowed with a  Lie superalgebra structure 
$[\ , \ ]: \End_\C(V)\times\End_\C(V)\lra \End_\C(V)$, 
which is a bilinear map defined by 
\[
[X, Y]= X Y - (-1)^{[X][Y]} Y X
\] 
for any homogeneous $X, Y\in \End_\C(V)$. Then $ \End_\C(V)$ together with $[\ ,\ ]$ is called the   {\em general linear Lie superalgebra} of $V$ and is denoted by $\gl(V)$.  

If $V$ is finite dimensional, there are non-negative integers $p$ and $q$ such that $\dim V_{\bar 0}=p$ and  $\dim V_{\bar 1}=q$. By choosing a homogeneous basis, we can identify $V$ with $\C^{p|q}=\C^p\oplus \C^q$.
We shall also  write   $\gl_{p|q}=\gl_{p|q}(\C)$ for $\gl(\C^{p|q})$.

The standard basis for $\C^{p|q}$ is given by $\{e_1,e_2,...,e_{p+q}\}$ where
\begin{equation}\label{eq_stbcpq}
e_1 = \begin{bmatrix}
           1 \\
           0\\
           0\\
           \vdots \\
           0
         \end{bmatrix},\quad
e_2 = \begin{bmatrix}
           0 \\
           1\\
           0\\
           \vdots \\
           0
         \end{bmatrix},        \quad 
\dots, \quad 
e_{p+q-1} = \begin{bmatrix}
           0 \\
           \vdots \\
           0\\
           1\\
           0
         \end{bmatrix}, \quad 
e_{p+q} = \begin{bmatrix}
           0 \\
           \vdots \\
           0\\
           0\\
           1
         \end{bmatrix}.
\end{equation}Then
$\{e_1, e_2, \dots, e_p\}$ is a basis for the even subspace  $\C^p$ and  $(e_{p+1}, e_{p+2}, \dots, e_{p+q})$ is a basis for the odd subspace $\C^q$. Denote by $E_{a b}$, for all $a, b=1, 2, \dots, p+q$, the matrix units relative to the standard basis. Then 
$
E_{a b} e_c = \delta_{b c} e_a$ for all $a, b, c. 
$

For positive integers $a$ and $b$, let $\rM_{ab}=\rM_{ab}(\C)$ be the space of all $a\times b$ complex matrices, $0_{ab}$ the $a\times b$ zero matrix and $\rM_a=\rM_{aa}$. Write $\fg= \gl_{p|q}$, and we shall identify each element $T\in\gl_{p|q}$ with the matrix which represents $T$ with respect to the standard basis of $\C^{p+q}$. So as vector spaces, we have $\fg=\rM_{p+q}$, 
\begin{align}\label{eq_g0}
\fg_{\bar 0} &=\left\{
\left(\begin{array}{l|l}
         A&0_{pq}\\ \hline
         0_{qp}&D
         \end{array}\right): A\in\rM_p,\ D\in\rM_q
\right\}, \\ 
\fg_{\bar 1} &=\left\{
\left(\begin{array}{l|l}
         0_{pp}& B\\ \hline
         C &0_{qq}
         \end{array}\right): B\in\rM_{pq},\ C\in\rM_{qp}
\right\}.   
\end{align}
The matrix units $E_{a b}$ form a homogeneous basis for $\fg$. 
We fix for $\fg$ the standard Borel subalgebra $\fb_{p|q}$ 
and Cartan subalgebras $\fh_{p|q}\subset \fb_{p|q}$, 
where $\fb_{p|q}$ consists of all upper triangular matrices in $\fg$, and $\fh_{p|q}$ consists of all diagonal matrices in $\fg$. We also let $\fu_{p|q}$ be the subalgebra of $\fg$ consisting of all strictly upper triangluar matrices in $\fg$. Then the Borel subalgebra $\fb_{p|q}$ can be written as a direct sum 
\[\fb_{p|q}=\fh_{p|q}\oplus\fu_{p|q}.\]

\begin{remark}
Note that  $\gl_{p|0}\cong \gl_p$, and $\fb_p:=\fb_{p|0}$ 
and $\fh_p:=\fh_{p|0}\subset \fb_{p|0}$ are the standard Borel and Cartan subalgebras of $\gl_p$ respectively. Moreover, $\fb_p=\fh_p\oplus \fu_p$ where $\fu_p:=\fu_{p|0}$. We
also have $\gl_{0|q}\cong \gl_q$.
\end{remark}

Highest weight modules for $\fg$ will be defined relative to the  Borel sublgebra $\fb_{p|q}$. 
Irreducible highest weight modules are uniquely characterized by their highest weights, which are elements of 
the dual space $\fh_{p|q}^*$  of $\fh_{p|q}$. The irreducible highest weight module with  highest weight
$\lambda\in \fh_{p|q}^*$ is finite dimensional if and only if $\lambda(E_{a a})- \lambda(E_{a+1,  a+1})\in \Z_+$ for all $a\ne p$.

We will denote an element $\lambda\in \fh_{p|q}^*$ as 
$\lambda=(\lambda_1, \lambda_2, \dots, \lambda_{p+q})$, where    
$\lambda=\lambda(E_{a a})$ for all $a=1, 2, \dots, p+q$.

Let us introduce a subset of $\fh_{p|q}^*$, which will play an important role in the rest of this paper.
Denote by $\Lambda^+$ the set of partitions, i.e., the set of sequences
$F=(f_1, f_2, \dots)$ with $f_i\in\Z_+$ for all $i$ such that $f_1\ge f_2\ge \dots$, and $f_j=0$ for all $j>>1$. 
We depict $F\in \Lam^+$ as a Young diagram with $f_i$ boxes in the $i$-th row for each $i$.
(Recall that a Young diagram is an array of square boxes
arranged in left-justified horizontal rows, with each row no longer
than the one above it \cite{Fu}. We will often abuse notation and write an element $F=(f_1,f_2,...)\in\Lambda^+$ as $F=(f_1,f_2,...,f_r)$ if $f_i=0$ for all $i>r$. For example, the element $(1,0,0,...)$ of $\Lambda^+$ will sometimes  be written as $(1)$, and it represents the Young diagram which has exactly $1$ box.) The depth of $F$, denoted by $\ell(F)$, is the number of rows in the Young diagram of $F$.  Let
\[\Lambda^+_n=\{D\in \Lambda^+: \ell(D)\le n\}.\] 
We also let
\[
\Lambda^+_{p|q}= \left\{F=(f_1, f_2, \dots)\in\Lambda^+ :   f_{p+1}\le q \right\}
\]
and call it the $(p, q)$-hook in $\Lambda^+$.
Note that $\Lambda^+_p\subset \Lambda^+_{p|q}$. 

Given $F=(f_1, f_2, \dots)\in\Lambda^+_{p|q}$, we write $f'_i$ for the length of the $i$-th column of the Young diagram of $F$ (in particular, $\ell(F)=f'_1$), and let $F^\sharp$ be the
 $(p+q)$-tuple defined by
\beq
F^\sharp= (  f^\sharp_{1}, f^\sharp_{2}, \dots, f^\sharp_{p+q}), 
\eeq
where for $1\le i\le p+q$,
\[f^\sharp_i=\left\{\begin{array}{ll}
f_i&1\le i\le p\\
\mathrm{max}(f'_{i-p} - p, 0)&p+1\le i\le p+q. 
\end{array}\right.\]  Let
 \[\Lambda_{p|q}^{+\sharp} =\{F^\sharp: F\in \Lambda^+_{p|q}\},\]  and we identify it with a subset of $\fh_{p|q}^*$ via the inclusion $\psi_{p|q}: \Lambda_{p|q}^{+\sharp}\lra  \fh_{p|q}^*$ defined by, for all $F=(f_1, f_2, \dots)\in\Lambda^+_{p|q}$, 
\beq\label{eq:wt-id}
\psi_{p|q}( F^\sharp)= \psi^{F^\sharp}_{p|q}, 
\eeq
where $\psi^{F^\sharp}_{p|q}(E_{ii})= f^\sharp_{i}$ for all $1\le i\le p+q$.

Now, for any $F\in {\Lambda^+_{p|q}}$, we denote  by $L_{p|q}^F$ the irreducible $\gl_{p|q}$-representation with highest weight 
$\psi^{F^\sharp}_{p|q}$. We call the representation of $\gl_{p|q}$ corresponding to such an $L_{p|q}^F$ an irreducible {\bf\em polynomial representation} in view of the well-known facts that tensor powers of $\C^{p|q}$ are semi-simple, and their simple submodules are precisely $L_{p|q}^F$ for all $F\in {\Lambda^+_{p|q}}$ (see, e.g., \cite{BR, GZ}).  Note in particular that $L_{p|q}^{(1)}=\C^{p|q}$.
\begin{remark}
The easiest way to see the semi-simplicity of the tensor powers of $\C^{p|q}$ is by observing the fact \cite{GZ} that $\C^{p|q}$ is unitarizable with respect to the compact real form of $\gl_{p|q}$. However, one should be warned that tensor products of polynomial representations with their duals  are not semi-simple in general.  
\end{remark}

\subsubsection{Realization of irreducible polynomial representations of $\gl_{p|q}$}
  Fix a positive integer $n$, and consider the natural action of 
$\gl_n\oplus\mathfrak{gl}_{p|q}$ on $\C^n\otimes \C^{p|q}$, where $\C^n$ is purely even. 
The action naturally extends to 
\[
\calr=\bigs(\C^n\otimes \C^{p|q}),
\]
the $\Z_2$-graded symmetric algebra on $\C^n\otimes \C^{p|q}$. 
We have the following result. 
\begin{theorem}[The $(\gl_n, \gl_{p|q})$-duality] \label{thm:duality-super}
As a $\gl_n\oplus\mathfrak{gl}_{p|q}$-module, $\calr$ has the following multiplicity free decomposition into irreducible represenations
\begin{equation}\label{eq_rdecom}
\calr\cong\bigoplus_{F\in\Lambda^{+}_{n,p|q}}\rho^F_n\otimes L^F_{p|q},
\end{equation}
where 
$
\Lambda^{+}_{n,p|q}=\Lambda^{+}_n\cap \Lambda^{+}_{p|q}, 
$
and  $\rho^F_n$ is the irreducible $\gl_n$ representation with highest weight determined by $F$ in the usual way (see Remark \ref{rmk:hw} below). 
\end{theorem}

\begin{remark}\label{rmk:hw}
Since any  $F=(f_1,f_2,...)\in \Lam^+_n$ satisfies $f_{n+i}=0$ for all $i\ge 1$, we define the map 
\[
\psi_n: \Lam^+_n\lra \fh_n^*, \quad F\mapsto \psi_n^F =  (f_1, f_2, \dots, f_n), 
\]
which is the special case $\psi_{n|0}$ of the map $\psi_{p|q}$ defined in equation \eqref{eq:wt-id}. 
Then  the highest weight of $\rho^F_n$ is equal to $ \psi_n^F$. 
 \end{remark}
 
Theorem \ref{thm:duality-super} was already known to Howe in the 70s \cite{H1}. 
More recent treatments of it can be found in \cite{Se} and \cite{CLZ, CW, CZ}. 
It reduces to $(\gl_n, \gl_p)$-duality when $q=0$:
\begin{equation}\label{eq_glnpd}
S(\C^n\otimes\C^p)\cong\bigoplus_{D\in\Lambda^+_{\min(n,p)}}\rho^D_n\otimes\rho^D_p,
\end{equation}
and to  skew $(\gl_n, \gl_q)$-duality when $p=0$:
\begin{equation}\label{eq_glnqsd}
\Lambda(\C^n\otimes\C^q)\cong\bigoplus_{E\in R_{n,q}}\rho^{E^t}_n\otimes\rho^E_q,
\end{equation}
where
\begin{equation}\label{eq_rnqdef}
R_{n,q}=\{F=(f_1,f_2,...)\in \Lambda^+: f_1\le n,\ f_{q+1}=0\},
\end{equation}
i.e., $R_{n,q}$ is  the set of all Young diagrams with at most $n$ columns and at most $q$ rows.
 (see \cite{H, H1}).

Next, denote by  $\fu_n$
 the subalgebra of $\gl_n$ consisting of all $n\times n$ strictly upper triangular matrices.  Let \[\calr^{\fu_n}=\{\xi\in\calr:\ X.\xi=0\ \forall X\in\fu_n\}.\]
 Then $\calr^{\fu_n}$ is a subalgebra of $\calr$, as for any $\zeta, \xi\in \calr$, we have 
 \[
 X.(\zeta\xi)= X.\zeta \xi + \zeta X.\xi =0, \quad  \forall  X\in\fu_n.
 \]
 It is  a module for $\fh_n\oplus\mathfrak{gl}_{p|q}$, and can be decomposed as
 \begin{equation}\label{eq_qdecom}
 \calr^{\fu_n}\cong\bigoplus_{F\in\Lambda^+_{n,p|q}}(\rho^F_n)^{\fu_n}\otimes L^F_{p|q}, 
 \end{equation}
where $(\rho^F_n)^{\fu_n} $ is the subspace of $\rho^F_n$ consisting of all vectors annihilated by the operators from $\fu_n$.
For $F\in \Lambda^+_{n,p|q} $, we have  $\dim (\rho^F_n)^{\fu_n}=1$ and  $h.v=\psi^F_n(h)v$    for all $v\in (\rho^F_n)^{\fu_n} $ and all  $h\in\fh_n$, i.e., the nonzero vectors in $(\rho^F_n)^{\fu_n} $ are
 $\gl_n$ highest weight vectors of weight   $\psi^F_n$. 
 Since $\dim (\rho^F_n)^{\fu_n}=1$, we have
 \begin{equation}\label{eq_lepq}
 (\rho^F_n)^{\fu_n}\otimes L^F_{p|q}\cong  L^F_{p|q}
 \end{equation}
 as modules for $\mathfrak{gl}_{p|q}$.
 Thus, we may realize the $\mathfrak{gl}_{p|q}$-module $L^F_{p|q}$ as the $\psi^F_n$-eigenspace of $\fh_n$ in $\calr^{\fu_n}$.

\subsection{Iterated Pieri rules for $\gl_n$}
In this subsection, we restrict ourselves to ordinary general linear Lie algebras and consider the cases  of Theorem \ref{thm:duality-super} with $q=0$ or $p=0$,  that is,  the $(\gl_n, \gl_p)$- duality or the skew $(\gl_n, \gl_q)$-duality. 
In particular,  
the special case with $p=1$ and $q=0$ yields the well-known fact that $\rho^{(d)}_n\cong S^d(\C^n)$,
the symmetric $d$-th power of the standard representation $\C^n$,  for any $d\in\Z_+$; 
and the case with $p=0$ and $q=1$ gives 
$\rho^{1^d}_n\cong \Lambda^d(\C^n)$, the exterior $d$-th power of $\C^n$, 
for any $d\le n$.

We shall state the {\em iterated Pieri rules} for $\gl_n$ which   play a crucial  role in the rest of this paper. Before we do that, we recall some definitions in combinatorics.
Let $D=(d_1,d_2,...)$ and $F=(f_1,f_2,...)$ be Young diagrams such that
 $D$ sits inside
$F$,  that is, $d_i\leq f_i$ for all $ i\geq 1$. By removing all boxes belonging to $D$, we obtain the
{\em skew diagram} $F/D$. (We allow $D$ to be the empty Young diagram with no box, i.e., $D=(0,0,...)$. In this case, we have $F/D=F$.) If we put a positive number in each box of
$F/D$, then it becomes a {\em skew tableau}. Let us denote this skew tableau by $T$. We say that the {\em
shape} of $T$ is $F/D$. If the entries of $T$
are taken from $[m]=\{1,2,...,m\}$, and $\alpha_j$ of them are $j$
for $1\leq j\leq m$, then we say the {\em content} of $T$ is $\alpha=(\alpha_1,...,\alpha_m)$. We also say that $T$
 is  {\em semistandard} if the
numbers in each row of  $T$  weakly increase from left-to-right, and
the numbers in each column of $T$ strictly increase from
top-to-bottom. 
The number of semistandard tableaux of shape $F/D$ and content $\alpha$ is denoted 
by $K_{F/D,\alpha}$ and is called a {\em Kostka number}(\cite{Fu}). (If $D$ is the empty Young diagram, then we write $K_{F,\alpha}$ for $K_{F/D,\alpha}$.) Finally,
 the {\em conjugate diagram}  of the Young diagram $D=(d_1,d_2,...)$ is the Young diagram $D^t$ having $d_j$ boxes in the $j$th column counting from left to right. 
 
\begin{proposition}[\textbf{Iterated Pieri rules}]\label{ispr}
Let $D,E\in\Lambda^{+}_n$.    
\begin{enumerate}
\item[(i)] If $\alpha_1,...,\alpha_p$ are non-negative integers, then 
\[
\dim \Hom_{\gl_n} \left( \rho^E_n, \rho^D_n \otimes  S ^{\alpha_1}(\C^n)
         \otimes \cdots \otimes
 S ^{\alpha_p}(\C^n) \right) = K_{E/D,\alpha}
\]
where $\alpha=(\alpha_1,...,\alpha_p)$.
\item[(ii)] If $\beta_1,...,\beta_q$ are non-negative integers not more than $n$, then
\[
\dim \Hom_{\gl_n} \left( \rho^E_n, \rho^D_n \otimes  \Lambda  ^{\beta_1}(\C^n)
         \otimes \cdots \otimes
 \Lambda  ^{\beta_q}(\C^n) \right) = K_{E^t/D^t,\beta},
\]
where $\beta=(\beta_1,...,\beta_q)$.
\end{enumerate}
 
\end{proposition}
\begin{proof} 
The proof of part (i) was given in \cite[ \S7.2]{L} and also in \cite[\S3.3]{HKL}, by combining the $(\gl_n, \gl_p)$-duality with the Pieri rule. 
Part (ii) was proven in \cite[Theorem 4]{KLW} in a similar spirit by combining the skew $(\gl_n, \gl_q)$-duality with the Pieri rule for skew symmetric tensors. 
\end{proof}


\section{Branching algebra for $(\gl_{p|q},\gl_{r|s}\oplus \fh_{r'}\oplus\fh_{s'})$  } \label{sec:problem1}
 
In this section, we shall construct a branching algebra for $(\gl_{p|q},\gl_{r|s}\oplus \fh_{r'}\oplus\fh_{s'})$ where $r+r'= p$ and $ s+s'= q$. By using the structure of this algebra, we obtain branching rules for  $(\gl_{p|q},\gl_{r|s}\oplus \fh_{r'}\oplus\fh_{s'})$ and also for 
$(\gl_{p|q},\gl_{r|s})$. By applying the result for the case $r=s=0$, we obtain the dimensions of the weight spaces in an irreducible polynomial representation of $\gl_{p|q}$.

  \subsection{Branching problems}\label{sec_glrs}
  From now on, we shall  denote $\gl_{p|q}$ by $\fg$. We first define two subalgebras of $\fg$ to be denoted by $\fg'$ and $\fm$ respectively.

  Let   $0\le r\le p$ and $0\le s\le q$ be non-negative integers, and let $\iota:\C^{r|s}\to\C^{p|q}$ be the injection defined by
\[(\begin{bmatrix}
           x_{1} \\
           x_{2} \\
           \vdots \\
           x_{r}
         \end{bmatrix},
         \begin{bmatrix}
           y_{1} \\
           y_{2} \\
           \vdots \\
           y_{s}
         \end{bmatrix})\to (\begin{bmatrix}
           x_{1} \\
           x_{2} \\
           \vdots \\
           x_{r}\\
           0\\
           \vdots\\
           0
         \end{bmatrix},
         \begin{bmatrix}
           y_{1} \\
           y_{2} \\
           \vdots \\
           y_{s}\\
           0\\
           \vdots\\
           0
         \end{bmatrix}).
         \]
 This induces an injection $\tilde{\iota}:\fgl_{r|s}\to\fg$ of Lie superalgebras defined as follows:
For any  $g=\left(\begin{array}{l|l}
         A&B\\ \hline
         C&D
         \end{array}\right)\in\fgl_{r|s}$,
 where the matrices $A,B,C,D$ are respectively $r\times r$, $r\times s$, $s\times r$ and $s\times s$, then
         \begin{equation}\label{eq_iotag}
         \tilde{\iota}(g)=\left(\begin{array}{c|c|c|c}
         A&0_{r,r'}&B&0_{r,s'}\\ \hline
         0_{r',r}&0_{r',r'}&0_{r',s}&0_{r',s'}\\ \hline
         C&0_{s,r'}&D&0_{s,s'}\\ \hline
         0_{s',r}&0_{s',r'}&0_{s',s}&0_{s',s'}
         \end{array}\right)\in\fg,\end{equation}
         where $0_{ab}$ denotes the $a\times b$ zero matrix, and 
         \begin{equation}
         r'=p-r,\quad s'=q-s.
         \end{equation}
Let
\begin{equation}\label{eq_gpdefn}
\fg'=\{ \tilde{\iota}(g): g\in\gl_{r|s}\},
\end{equation} which is a 
subalgebra of $\fg$ and $\fg'\cong\fgl_{r|s}$.
Let
\[\fh_{\fg'}=\tilde{\iota}(\fh_{r|s}),\quad \fu_{\fg'}=\tilde{\iota}(\fu_{r|s}),\quad\mbox{and}\quad\fb_{\fg'}=\tilde{\iota}(\fb_{r|s}).\]
Then $\fh_{\fg'}$ is a Cartan subalgebra, $\fb_{\fg'}$ is a Borel subalgebra of $\fg'$
and $\fb_{\fg'}=\fh_{\fg'}\oplus\fu_{\fg'}$.

We also let $\fm$ be the subalgebra of $\fg$ consisting of all matrices  of the form \eqref{eq_iotag}
 \begin{equation}\label{eq_mdefine}
         \left(\begin{array}{c|c|c|c}
         A&0_{r,r'}&B&0_{r,s'}\\ \hline
         0_{r',r}&H_{r'}&0_{r',s}&0_{r',s'}\\ \hline
         C&0_{s,r'}&D&0_{s,s'}\\ \hline
         0_{s',r}&0_{s',r'}&0_{s',s}&H_{s'}
         \end{array}\right)\end{equation}
         where all the submatrices are as in equation \eqref{eq_iotag} except that $H_{r'}$ and $H_{s'}$ are diagonal matrices of sizes $r'\times r'$ and $s'\times s'$  respectively. Then it is clear that 
\[
\fm\cong \fg'\oplus \fkh_{r'}\oplus \fkh_{s'}\cong \gl_{r|s}\oplus \fkh_{r'}\oplus \fkh_{s'}.
\]  
Let
 \[\fb_\fm=\fm\cap\fb_{p|q},\quad\fh_\fm=\fm\cap\fh_{p|q},\quad\fu_\fm=\fm\cap\fu_{p|q}.\]
 Then $\fb_m $ is a Borel subalgebra of $\fm$, $\fh_\fm$ is a Cartan subalgebra of $\fm$ and $\fb_\fm=\fh_\fm\oplus\fu_\fm$. Note that  
 \[\fh_{\fm} \cong \fh_{\fg'}\oplus\fh_{r'}\oplus \fh_{s'}
\cong \fh_{r|s}\oplus\fh_{r'}\oplus \fh_{s'},\quad\mbox{and}\quad
\fu_m=\fu_{\fg'}.\]

 In this section, we consider the following branching problems: 
 
\bigskip\noindent 
 {\bf
 Branching Problem:}    {\em For $F\in\Lambda^+_{p|q}$, 
 \begin{enumerate}
 \item[(i)] determine a decomposition of $L^F_{p|q}$ into irreducible $\fm$ modules, and
 \item[(ii)] determine a decomposition of $L^F_{p|q}$ into irreducible $\fg'$ modules.
  \end{enumerate}}

 \subsection{The branching algebra $\calq(\fg,\fm)$}\label{sec_bagm23}
 Let $n$ be a positive integer and consider the
 $\gl_n\oplus\fg$ module 
$\calr=\bigs(\C^n\otimes \C^{p|q})$.
  Let 
 \begin{equation}\label{eq_calqdef}
 \calq(\fg,\fm)=\calr^{\fu_n\oplus\fu_{\fm}}
 \end{equation}
 be the subalgebra of $\calr$ consisting of vectors annihilated by operators from $\fu_n\oplus\fu_{\fm}$.
 It is a module for $\fh_n\oplus\fh_\fm$.   
 Since  $\fh_n\oplus\fh_\fm$ is isomorphic $\Sigma\fh:=\fh_n\oplus\fh_{r|s}\oplus\fkh_{r'}\oplus\fkh_{s'}$, we shall instead regard  $\calq(\fg,\fm)$ as a $\Sigma\fh$ module and 
  describe the corresponding isotypic decomposition of $\calq(\fg,\fm)$. 
 
 \medskip
 For convenience, we introduce the following notation:
 \begin{notation} Let $U_1,U_2,...,U_k$ be complex vector spaces
 and  $U=U_1\oplus U_2\oplus\cdots\oplus U_k$. For each $1\le i\le u$, let $\psi_i:U_i\to\C$ be a linear functional. Then $ (\psi_1,\psi_2,...,\psi_k)$ shall denote the linear functional on $U$ defined by, for all $u=(u_1,u_2,...,u_k)\in U$,
 \begin{equation}\label{notation_psi}
(\psi_1,\psi_2,...,\psi_k)(u)=\psi_1(u_1)+\psi_2(u_2)+\cdots+\psi_k(u_k).
 \end{equation}
 \end{notation}

 We shall assume that $(r,s)\ne (0,0)$ here,   the case $(r,s)=(0,0)$ will be discussed in \S \ref{sec_wtm}.
Let
\begin{equation}\label{eq_omegadefn}
\Omega(\fg,\fm):=\Lambda^{+}_{n,p|q}\times\Lambda^+_{n,r|s}\times \Z^{r'}_{+}\times \Z^{s'}_{+}.
\end{equation}
For each $(F,D,\alpha,\beta)\in\Omega(\fg,\fm)$, let 
\[\calq(\fg,\fm)_{(F,D,\alpha,\beta)}=\left\{v\in \calq(\fg,\fm):h.v=(\psi^F_n,\psi^{D^\sharp}_{p|q},\psi^\alpha_{r'},\psi^\beta_{s'})(h)\ \forall h\in \Sigma\fh\right\},\] 
that is, $\calq(\fg,\fm)_{(F,D,\alpha,\beta)}$ is the $(\psi^F_n,\psi^{D^\sharp}_{p|q},\psi^\alpha_{r'},\psi^\beta_{s'})$-isotypic component for  $\Sigma\fh$ in $\calq(\fg,\fm)$.

Let us denote by $\C_{\psi^\alpha_{r'}}$ 
the $1$-dimensional $\fh_{r'}$-module with weight $\psi^\alpha_{r'}$, 
and similarly denote  by $\C_{\psi^\beta_{s'}}$ the $1$-dimensional $\fh_{s'}$-module with weight  $\psi^\beta_{s'}$. Then we have the following result:
\begin{proposition}\label{prop_qgmdecom}
The algebra $\calq(\fg,\fm)$ has a direct sum decomposition given by 
\begin{equation}\label{eq_qgmisodecom}
\calq(\fg,\fm)=\bigoplus_{(F,D,\alpha,\beta)\in\Omega(\fg,\fm)}\calq(\fg,\fm)_{(F,D,\alpha,\beta)}.
\end{equation}
Moreover, for each $(F,D,\alpha,\beta)\in\Omega(\fg,\fm)$, 
\begin{equation}\dim \calq(\fg,\fm)_{(F,D,\alpha,\beta)}
 =\dim\Hom_\fm(L^D_{r|s}\otimes\C_{\psi^\alpha_{r'}}\otimes\C_{\psi^\beta_{s'}}, L^F_{p|q}).\label{eq_dimqgm1}
 \end{equation}
\end{proposition}
\begin{proof}
Since all vectors in $\calq(\fg,\fm)$ are anihilated by the operators from $\fu_n$, any $\fh_n$ eigenvector in $\calq(\fg,\fm)$  is a $\gl_n$ highest weight vector, so that the corresponding eigencharacter is necessarily a dominant weight for $\gl_n$ and is of the form $\psi^F_n$ for some $F\in\Lambda^+_{n,p|q}$. For a similar reason, the eigencharacter of any eigenvector of $\fh_{r|s}$ in $\calq(\fg,\fm)$ is of the form $\psi^{D^\sharp}_{r|s}$ for some $D\in\Lambda^+_{n,r|s}$. Consequently, the direct sum \eqref{eq_qgmisodecom} 
is the isotypic decomposition of $\calq(\fg,\fm)$ with respect to the action by $\Sigma\fh$.

To prove \eqref{eq_dimqgm1},  we shall use equation \eqref{eq_rdecom} to obtain another  decomposition of $\calq(\fg,\fm)$.  By extracting vectors in $\calr$ which are anihilated by $\fu_n\oplus\fu_m\cong\fu_n\oplus\fu_{\fg'}$,  we obtain from equation \eqref{eq_rdecom} the direct sum
\begin{equation}\label{eq_qgmdecom1}
\calq(\fg,\fm)=\calr^{\fu_n\oplus\fu_{\fm}}\cong\bigoplus_{F\in\Lambda^{+}_{n,p|q}}(\rho^F_n)^{\fu_n}\otimes (L^F_{p|q})^{\fu_{\fg'}},
\end{equation}
where 
  $(L^F_{p|q})^{\fu_{\fg'}}$ is the space of vectors in $L^F_{p|q}$ annihilated by operators from $\fu_{\fg'}$, 
  which is a module for $\fh_{r|s}\oplus\fh_{r'}\oplus\fh_{s'}$.  
  
 We now fix $F\in\Lambda^+_{n,p|q}$. Then $(L^F_{p|q})^{\fu_{\fg'}}$ can be decomposed as
\begin{equation}\label{eq_lfpqugpdecom}
(L^F_{p|q})^{\fu_{\fg'}}= \bigoplus_{(D,\alpha,\beta)}(L^F_{p|q})^{\fu_{\fg'}}_{(D,\alpha,\beta)},
\end{equation}
where the direct sum is taken over all $(D,\alpha,\beta)\in \Lambda^+_{n,r|s}\times \Z^{r'}_{+}\times \Z^{s'}_{+}$ and
\[(L^F_{p|q})^{\fu_{\fg'}}_{(D,\alpha,\beta)}=\left\{v\in (L^F_{p|q})^{\fu_{\fg'}}: h.v=(\psi^{D^\sharp}_p, \psi^\alpha_{r'}, \psi^\beta_{s'})(h)v \forall h\in\fh_{r|s}\oplus\fh_{r'}\oplus\fh_{s'}\right\},\]
 i.e.,  
 the space of all $\fm$ highest weight vectors in $L^F_{p|q}$ of weight $(\psi^{D^\sharp}_n, \psi^\alpha_{r'},\psi^\beta_{s'})$. Consequently,
 \begin{equation}\label{eq_dimlfpqfgp}
 \dim (L^F_{p|q})^{\fu_{\fg'}}_{(D,\alpha,\beta)}=\dim \Hom_\fm(L^D_{r|s}\otimes\C_{\psi^\alpha_{r'}}\otimes\C_{\psi^\beta_{s'}}, L^F_{p|q}),
 \end{equation}
 which is the multiplicity of $L^D_{r|s}\otimes\C_{\psi^\alpha_{r'}}\otimes\C_{\psi^\beta_{s'}}$ in $L^F_{p|q}$.

 Using equations \eqref{eq_qgmdecom1} and \eqref{eq_lfpqugpdecom}, we obtain 
 \begin{align}
\calq(\fg,\fm) 
&\cong\bigoplus_{(F,D,\alpha,\beta)\in\Omega(\fg,\fm)}(\rho^F_n)^{\fu_n}\otimes (L^F_{p|q})^{\fu_{\fg'}}_{(D,\alpha,\beta)
}.\label{eq_qgmdecom2}\end{align}

 By comparing equations \eqref{eq_qgmisodecom} and \eqref{eq_qgmdecom2}, we see that for each $(F,D,\alpha,\beta)\in\Omega(\fg,\fm)$,
\begin{equation}\label{eq_qgmcong1}
\calq(\fg,\fm)_{(F,D,\alpha,\beta)}\cong (\rho^F_n)^{\fu_n}\otimes (L^F_{p|q})^{\fu_{\fg'}}_{(D,\alpha,\beta)}\end{equation}
as a module for $\Sigma\fh$.  
 Since $\dim (\rho^F_n)^{\fu_n}=1$ and by equation \eqref{eq_dimlfpqfgp}, we obtain
\[
\dim \calq(\fg,\fm)_{(F,D,\alpha,\beta)} 
=\dim   (L^F_{p|q})^{\fu_{\fg'}}_{(D,\alpha,\beta)} 
=\dim\Hom_\fm(L^D_{r|s}\otimes\C_{\psi^\alpha_{r'}}\otimes\C_{\psi^\beta_{s'}}, L^F_{p|q}).\]
This proves \eqref{eq_dimqgm1}, completing the proof of the proposition. 
 \end{proof}

 \begin{remark}\label{rm_ba1}
 We see from Proposition \eqref{prop_qgmdecom}  that part of the branching rule for $(\fg,\fm)$ can be deduced from  the structure of  the algebra $\calq(\fg,\fm) $. In particular, the dimension of  each  of the subspaces  $\calq(\fg,\fm)_{(F,D,\alpha,\beta)}$ of $\calq(\fg,\fm) $   is the multiplicity 
 of an irreducible $\fm$ representation in an irreducible $\fg$ representation. 
 In view of this property of $\calq(\fg,\fm)$, we call $\calq(\fg,\fm) $ a {\bf\em branching algebra} for $(\fg,\fm)$.
 \end{remark}

\subsection{Branching rules for $(\fg,\fm)$ and $(\fg,\fg')$}\label{sec_brgm23} We now  express $\dim \calq(\fg,\fm)_{(F,D,\alpha,\beta)}$ in terms of Kostka numbers and skew Kostka numbers.  

\begin{proposition}\label{prop_qgm2}
For any $(F,D,\alpha,\beta)\in\Omega(\fg,\fm)$,  let $N_{(F,D,\alpha,\beta)}$ be the non-negative integer defined by
\begin{equation}\label{eq_nfdabdefn}
N_{(F,D,\alpha,\beta)}=\sum_EK_{E/D,\alpha}K_{F^t/E^t,\beta}
\end{equation}
where the sum is taken over all Young diagrams $E$ (note that $K_{E/D,\alpha}K_{F^t/E^t,\beta}\ne 0 $ only for finitely many $E$). Then  
\begin{equation} 
\dim \calq(\fg,\fm)_{(F,D,\alpha,\beta)}=N_{(F,D,\alpha,\beta)}.\label{eq_dimqgm2}
 \end{equation}
\end{proposition}
\begin{proof} 
Write $\C^p=\C^r\oplus\C^{r'}$ and $\C^q=\C^s\oplus\C^{s'}$. Then
$\C^{p|q}= \C^r\oplus \C^{r'}\oplus\C^{s}\oplus\C^{s'} \cong \C^{r|s}\oplus\C^{r'|s'}$, 
and hence
$\CR=\BS(\C^n\ot\C^{p|q}) \cong \BS(\C^n\ot\C^{r|s}\oplus\C^n\ot\C^{r'|s'})$. 
Recall from Section \ref{sect:super-facts} the functorial property of the supersymmetric algebra 
that $\BS(V\oplus V')=\BS(V)\otimes \BS(V')$ for any $\Z_2$-graded vector spaces $V$ and $V'$. 
This immediately leads to  
\[
\baln
\CR&\cong \BS(\C^n\ot\C^{r|s})\ot \BS(\C^n\ot\C^{r'|s'}).
\ealn
\]
Since $\bigs(\C^n\otimes\C^{r'|s'})\cong S(\C^n\otimes\C^{r'})\otimes\Lambda(\C^n\otimes\C^{s'})$ (see Section \ref{sect:super-facts}), we obtain
\begin{align}\label{eq:interm}
\calr  &\cong \bigs(\C^n\otimes\C^{r|s})\otimes S(\C^n\ot\C^{r'})
 \otimes \Lambda(\C^n\ot\C^{s'}). 
 \end{align}
This is an isomorphism of modules for $\gl_n\times \fm$, whose direct summands 
$\gl_n$, $\gl_{r|s}$, $\fh_{r'}$ and $\fh_{s'}$ act on both sides in the obvious way.

Next, we have  
  $
  S(\C^n\ot\C^{r'})\cong S(\C^n)^{\otimes r'}
  $  
  and 
  $\Lambda(\C^n\ot\C^{s'})
 \cong\Lambda(\C^n)^{\otimes s'}. 
 $
By using the $(\gl_n, \gl_1)$-duality and its skew version, we obtain
 \[
 S(\C^n)\cong\bigoplus_{a\in\Z_{+}} S^a(\C^n)\ot\rho_1^{(a)}, \quad 
  \Lambda(\C^n)\cong\bigoplus_{b=0}^n   \Lambda ^{b}(\C^n)
                           \otimes\rho^{(b)}_1. 
 \] 
  Thus we have the following isomorphisms of modules for $\gl_n\times\fh_{r'}$ and  $\gl_n\times\fh_{s'}$ respectively:
 \begin{align}
 &S(\C^n\otimes \C^{r'})\cong \bigotimes_{i=1}^{r' }\left(\bigoplus_{\alpha_i\in\Z_{+}}   S ^{\alpha_i}(\C^n)
                           \otimes\rho^{(\alpha_i)}_1\right) 
\cong \bigoplus_{\alpha\in\Z_{+}^{r'}}   \left(\bigotimes_{i=1}^{r' } S ^{\alpha_i}(\C^n)\right)
                           \otimes\C_{\psi_{r'}^\alpha},\label{eq_scncrp} \\
&\Lambda(\C^n\otimes \C^{s'})\cong   \bigotimes_{j=1}^{s'} \left(\bigoplus_{\beta_j=0}^n   \Lambda ^{\beta_j}(\C^n)
                           \otimes\rho^{(\beta_j)}_1\right) 
\cong \bigoplus_{\beta\in\Z_{+}^{s'}} \left(\bigotimes_{j=1}^{s'} \Lambda ^{\beta_j}(\C^n) \right)
                           \otimes\C_{\psi_{s'}^{\beta}},\label{eq_lcncsp}
\end{align}  
where $\alpha=(\alpha_1, \alpha_2, \dots, \alpha_{r'})$, $\beta=(\beta_1, \beta_2, \dots, \beta_{s'})$, 
and $\Lambda ^{\beta_j}(\C^n) =0$ if $\beta_j>n$ for any $j$. Moreover, 
$ \C_{\psi^\alpha_{r'}}$ and $\C_{\psi^\alpha_{s'}}$ are one-dimensional modules 
of  $\fkh_{r'}$ and $\fkh_{s'}$ respectively, satisfying 
  \begin{align*}
  t_1.f_1&=\psi^\alpha_{r'}(t_1)f_1 ,\qquad
  t_2.f_2=\psi^\beta_{s'}(t_2)f_2 ,
  \end{align*} 
  for all $t_1 \in\fkh_{r'}$, $f_1\in\C_{\psi^\alpha_{r'}}$,
  $t_2 \in\fkh_{s'}$ and $f_2\in\C_{\psi^\beta_{s'}}$.

It now follows from equations \eqref{eq:interm}, \eqref{eq_rdecom}, \eqref{eq_scncrp} and \eqref{eq_lcncsp} that  
 \begin{align}
 \calr &=\bigoplus_{(D,\alpha,\beta)\in\Lambda^+_{n,r|s}\times\Z^{r'}_{+}\times \Z^{s'}_{+}}  \calr_{(D,\alpha,\beta)}\label{eq_grading}
\end{align} 
where for each $(D,\alpha,\beta)\in\Lambda^+_{n,r|s}\times\Z^{r'}_{+}\times\Z^{s'}_{+}$, 
\[\calr_{(D,\alpha,\beta)}\cong T(D,\alpha,\beta)\otimes L^D_{r|s}\otimes \C_{\psi^\alpha_{r'}}\otimes\C_{\psi^\beta_{s'}},\]
and $T(D,\alpha,\beta)$ is the $\gl_n$ module defined by
\[T(D,\alpha,\beta)= 
\rho^D_n\otimes\left(\bigotimes_{i=1}^{r'}   S ^{\alpha_i}(\C^n)\right)
\otimes   \left(\bigotimes_{j=1}^{s'}  \Lambda ^{\beta_j}(\C^n)\right).\]
 (We agree that $T(D,\alpha,\beta)=0$ if $j>n$ for some $1\le j\le s'$.) 
   
Now $T(D,\alpha,\beta)$ can be analyzed using the Iterated Peiri rules. 
By first applying part (i) and then part (ii) of Proposition \ref{ispr},  we obatin 
\begin{align*}
T(D,\alpha,\beta)
&\cong
\left(\bigoplus_EK_{E/D,\alpha}\rho^E_n\right)
\otimes   \left(\bigotimes_{j=1}^{s'}  \Lambda ^{\beta_j}(\C^n)\right)\\
&\cong
\bigoplus_EK_{E/D,\alpha}\left(\bigoplus_F K_{F^t/E^t,\beta}\rho^F_n \right)\\
&\cong
\bigoplus_F\left(\bigoplus_EK_{E/D,\alpha}K_{F^t/E^t,\beta}\right) \rho^F_n \\
&\cong
\bigoplus_F N_{(F,D,\alpha,\beta)} \rho^F_n.
\end{align*}
By extracting the vectors in  $T(D,\alpha,\beta)$ which are anihilated by $\fu_n$, we obtain
\begin{equation}\label{eq_tdabdecom}
T(D,\alpha,\beta)^{\fu_n} \cong
\bigoplus_FN_{(F,D,\alpha,\beta)} (\rho^F_n)^{\fu_n}.\end{equation}
We now fix $F\in\Lambda^+_{n,p|q}$ and let $T(D,\alpha,\beta)^{\fu_n}_F$ be the space of vectors $v$ in 
$T(D,\alpha,\beta)^{\fu_n}$ for which $a.v=\psi^F_n(a)v$ for all $a\in\fh_n$. Then since $\dim (\rho^F_n)^{\fu_n}=1$ and  by equation \eqref{eq_tdabdecom},
\begin{equation}\label{eq_dimtdab}\dim T(D,\alpha,\beta)^{\fu_n}_F=N_{(F,D,\alpha,\beta)}.\end{equation}
Moreover, 
\[\calq_{(F,D,\alpha,\beta)}\cong T(D,\alpha,\beta)^{\fu_n}_F\otimes(L^D_{r|s})^{\fu_{r|s}}\otimes \C_{\psi^\alpha_{r'}}\otimes\C_{\psi^\beta_{s'}}\]
as a module for $\fh_n\oplus\fh_{r|s}\oplus\fh_{r'}\oplus\fh_{s'}$. It follows from this and  equation   \eqref{eq_dimtdab} that 
\begin{align*}
\dim\calq_{(F,D,\alpha,\beta)}&=\dim  \left(T(D,\alpha,\beta)^{\fu_n}_F\otimes(L^D_{r|s})^{\fu_{r|s}}\otimes \C_{\psi^\alpha_{r'}}\otimes\C_{\psi^\beta_{s'}}\right)\\
&=\dim T(D,\alpha,\beta)^{\fu_n}_F 
\dim(L^D_{r|s})^{\fu_{r|s}}
\dim  \C_{\psi^\alpha_{r'}}  \dim\C_{\psi^\beta_{s'}}  \\
&=\dim T(D,\alpha,\beta)^{\fu_n}_F \\
&=N_{(F,D,\alpha,\beta)}.
\end{align*}
This completes the proof. 
\end{proof}

 \begin{remark}\label{rm_ba2} We have seen earlier that the structure of $\calq(\fg,\fm)$ encodes a branching rule for $(\fg,\fm)$. The proof of Proposition \eqref{prop_qgm2} shows that the dimension of 
 $\calq_{(F,D,\alpha,\beta)}$ is equal to the multiplicity of $\rho^F_n$ in the $k$-fold tensor product $T(D,\alpha,\beta)$ of $\gl_n$ representations, where $k=r'+s'+1$. So the structure of  $\calq(\fg,\fm) $ also encodes information on how such tensor products decomposes into irreducible $\gl_n$ representations, which can be viewed as a branching rule from $\gl_n^{k}$ to its diagonal subalgebra  $\Delta^{(k-1)}(\gl_n)$ (see Remark \ref{rmk:co-mult} below for explanation).
  Hence, the structure of  $\calq(\fg,\fm)$  encodes two sets of branching rules connected by a reciprocity law. 
 Following \cite{HTW1,HTW2}, we shall call $\calq(\fg,\fm)$
 a  {\bf\em reciprocity algebra}.
 \end{remark}
 
\begin{remark}\label{rmk:co-mult}
%
 Set $\Delta^{(k-1)}(\gl_n)=\left\{\Delta^{(k-1)}(X)=(\underbrace{X, X, \dots, X}_k):  X\in\gl_n\right\}.$
When $\gl_n$ is regarded as embedded in its universal enveloping algebra in the canonical way, one always writes
$\Delta^{(k-1)}(X)=\sum_{i=0}^{k-1} \underbrace{1\ot\dots\ot 1}_i\ot  X\ot \underbrace{1\ot\dots\ot 1}_{k-1-i}$.
 \end{remark}

We have the following result
 
 \begin{theorem}\label{thm:brg-m} Recall that $r'=p-r$ and $s'=q-s$.   
 Let $F\in\Lambda^+_{n,p|q}$.
 
   \begin{enumerate}
  \item[(i)]  {\rm\bf(Branching from $\gl_{p|q}$ to $\gl_{r|s}\oplus \fh_{r'}\oplus\fh_{s'})$}As a representation of $\fm$,
    \[
  L^F_{p|q}=\bigoplus_{(D,\alpha,\beta)\in\Lam^+_{n, r|s}\times \Z^{r'}_{+}\times \Z^{s'}_{+}}N_{(F,D,\alpha,\beta)}\ L^D_{r|s}\otimes \C_{\psi^\alpha_{r'}}\otimes\C_{\psi^\beta_{s'}}.
  \]
  
  \item[(ii)] {\rm\bf(Branching from $\gl_{p|q}$ to $\gl_{r|s} $}
  As a representation of $\fg'$,
    \[
  L^F_{p|q}=\bigoplus_{D\in\Lam^+_{n, r|s} }\widetilde{N}_{(F,D)}\ L^D_{r|s} .
  \]
  where for each $D\in\Lam^+_{n, r|s} $, 
   \begin{equation}\label{eq_tnfd}
   \widetilde{N}_{(F,D)}=\sum_{(\alpha,\beta)\in\Z^{r'}_{+}\times \Z^{s'}_{+}}N_{(F,D,\alpha,\beta)}.
   \end{equation}
  \end{enumerate}
  
  \end{theorem}
\begin{proof}   
 This follows from equations \eqref{eq_dimqgm1} and \eqref{eq_dimqgm2}.
 \end{proof}

\subsection{Weight multiplicities of $L_{p|q}^F$}\label{sec_wtm}
We now consider the case  $r=s=0$.  In this case, we have 
\[
\fg'=0,\quad \fm=\fh_\fm=\fh_p\oplus\fh_q,\quad
\calq(\fg,\fm)=\calr^{\fu_n}, \quad
\Omega(\fg,\fm)=\Lambda^+_{n,p|q}\times\Z^p_{+}\times\Z^q_{+},
\]
and $\calq(\fg,\fm)$ is a module for $\Sigma\fh:=\fh_n\oplus\fh_p\oplus\fh_q$.
Note that $\Omega(\fg,\fm)$ is the direct product of only $3$ sets instead of $4$ since $\Lambda^+_{n,r|s}=\emptyset$, and 
the algebra $\calq(\fg,\fm)$ decomposes into the direct sum 
\[\calq(\fg,\fm)=\bigoplus_{(F,\alpha,\beta)\in\Omega(\fg,\fm)}\calq(\fg,\fm)_{(F,\alpha,\beta)},\]
where for each $(F,\alpha,\beta)\in\Omega(\fg,\fm)$, 
\[\calq(\fg,\fm)_{(F,\alpha,\beta)}=\left\{v\in \calq(\fg,\fm):h.v=(\psi^F_n,\psi^\alpha_p,\psi^\beta_q)(h)\ \forall h\in \Sigma\fh\right\}.\] 
Except these minor differences, the arguments in \S \ref{sec_bagm23} and \ref{sec_brgm23} remain valid and lead to the following results:

 \begin{corollary} \label{cor:CZ-basis}
 Let $F\in\Lambda^+_{n,p|q}$.
\begin{enumerate}
\item[(i)] For $(\alpha,\beta)\in\Z^p_{+}\times\Z^q_{+}$, 
let 
\[(L^F_{p|q})_{(\alpha,\beta)}=\left\{v\in L^F_{p|q}:\ (a,b).v=(\psi^\alpha_p(a)+\psi^\beta_q(b))v\ \forall (a,b)\in \fh_p\oplus\fh_q\right\},\]
the $(\psi^\alpha_p,\psi^\beta_q)$-weight space of $L^F_{p|q}$.
Then  
we have
\[\calq(\fg,\fm)_{(F,\alpha,\beta)}\cong (\rho^F_n)^{\fu_n}\otimes (L^F_{p|q})_{(\alpha,\beta)}\]
as a module for $\Sigma\fh$, and 
\begin{equation}\label{eq_wtmutlfpq}
\dim (L^F_{p|q})_{(\alpha,\beta)}=N'_{(F,\alpha,\beta)}
\end{equation}
where 
 \beq
N'_{(F,\alpha,\beta)}=\sum_EK_{E,\alpha}K_{F^t/E^t,\beta} \label{eq:mult-formula}
\eeq
where the sum is taken over all Young diagrams $E$.

\item[(ii)]  The dimension of $L^F_{p|q}$ is given by
\begin{equation}\label{eq_dimlfpq}
\dim L^F_{p|q}=\widetilde{N}^\prime_F,
\end{equation}
where 
\[\widetilde{N}'_F= \sum_{(\alpha,\beta)\in\Z^p_{+}\times\Z^q_{+}}N'_{(F,\alpha,\beta)}.\]
\end{enumerate}
 
 \end{corollary}
 
\begin{remark} 
Part (i) of Corollary \ref{cor:CZ-basis} gives the weight multiplicities of any  irreducible polynomial $\gl_{p|q}$-representation. We believe that the formula \eqref{eq_wtmutlfpq}, which expresses weight multiplicities in terms of Kostka numbers, is new. 
\end{remark}

\section{Branching algebra for $(\gl_{p|q},\gl_{r|s}\oplus\gl_{p-r|q-s})$} \label{glpglq}
 
 In this section, we shall construct a branching algebra for $(\gl_{p|q},\gl_{r|s}\oplus\gl_{r'|s'})$ where
$r+r'=p,\ s+s'=q$. By using the structure of this algebra, we deduce a branching rule for  $(\gl_{p|q},\gl_{r|s}\oplus\gl_{r'|s'})$.    In  particular, in the case $(r,s,r',s')=(p,0,0,q)$, we obtain a branching rule from $\gl_{p|q}$ to its even subspace $\gl_p\oplus\gl_q $.

  \subsection{Branching problem}\label{sec_bp2}

  We continue to denote $\fg=\gl_{p|q}$. Recall the subalgebra $\fg'$ of $\fg$ defined in equation \eqref{eq_gpdefn}. We now define another Lie supalgebra $\fgpp$ of $\fg$ as follows:   Let $\kappa:\C^{r'|s'}\to\C^{p|q}$ be the injection defined by
\[(\begin{bmatrix}
           x_{1} \\
           x_{2} \\
           \vdots \\
           x_{r'}
         \end{bmatrix},
         \begin{bmatrix}
           y_{1} \\
           y_{2} \\
           \vdots \\
           y_{s'}
         \end{bmatrix})\to (\begin{bmatrix}
         0\\
           \vdots\\
           0\\
           x_{1} \\
           x_{2} \\
           \vdots \\
           x_{r'} 
            \end{bmatrix},
         \begin{bmatrix}
         0\\
           \vdots\\
           0\\
           y_{1} \\
           y_{2} \\
           \vdots \\
           y_{s'} 
         \end{bmatrix}).
         \]
 This induces an injection $\tilde{\kappa}:\fgl_{r'|s'}\to\fg$ of Lie superalgebras defined as follows:
For any  $g=\left(\begin{array}{l|l}
         E&F\\ \hline
         G&H
         \end{array}\right)\in\fgl_{r'|s'}$,
 where the matrices $E,F,G,H$ are respectively $r'\times r'$, $r'\times s'$, $s'\times r'$ and $s'\times s'$, then
         \begin{equation}\label{eq_kappag}
         \tilde{\kappa}(g)=\left(\begin{array}{c|c|c|c}
         0_{r,r}&0_{r,r'}&0_{r,s}&0_{r,s'}\\ \hline
         0_{r',r}&E&0_{r',s}&F\\ \hline
         0_{s,r}&0_{s,r'}&0_{s,s}&0_{s,s'}\\ \hline
         0_{s',r}&G&0_{s',s}&H
         \end{array}\right)\in\fg.\end{equation}
         Let
\begin{equation}\label{eq:gp-defn}
\fgpp=\{ \tilde{\kappa}(g): g\in\gl_{r'|s'}\},
\end{equation} which is a 
subalgebra of $\fg$ and $\fgpp\cong\fgl_{r'|s'}$.
Let
\[\fh_{\fg''}=\tilde{\kappa}(\fh_{r'|s'}),\quad \fu_{\fg''}=\tilde{\kappa}(\fu_{r'|s'}),\quad\mbox{and}\quad\fb_{\fg''}=\tilde{\kappa}(\fb_{r'|s'}).\]
Then $\fh_{\fg''}$ is a Cartan subalgebra and $\fb_{\fg''}$ is a Borel subalgebra of $\fg''$, 
and $\fb_{\fg''}=\fh_{\fg''}\oplus\fu_{\fg''}$.

Next we let $\fq$ be the subspace of $\fg$ spanned by $\fg'\cup\fgpp$. Explicitly, $\fq$ consists of all matrices of the form
  \begin{equation}\label{eq_matfl}
         \left(\begin{array}{c|c|c|c}
         A&0_{r,r'}&B&0_{r,s'}\\ \hline
         0_{r',r}&E&0_{r',s}&F\\ \hline
         C&0_{s,r'}&D&0_{s,s'}\\ \hline
         0_{s',r}&G&0_{s',s}&H
         \end{array}\right)\in\fg.\end{equation}
 where all the submatrices are as in equations \eqref{eq_iotag}  and \eqref{eq_kappag}.
Then $\fq$ is a subalgebra of $\fg$, and 
\begin{equation}\label{eq_qcongglrsglrpsp}\fq\cong \fg'\oplus\fg''\cong\gl_{r|s}\oplus\gl_{r'|s'}.\end{equation}
Let
 \[\fb_\fq=\fq\cap\fb_{p|q},\quad\fh_\fq=\fq\cap\fh_{p|q},\quad\fu_\fq=\fq\cap\fu_{p|q}.\]
 Then $\fb_q$ is a Borel subalgebra of $\fq$, $\fh_\fq$ is a Cartan subalgebra of $\fq$ and $\fb_\fq=\fh_\fq\oplus\fu_\fq$. Moreover, we have
 \[\fh_{\fq} \cong \fh_{\fg'}\oplus\fh_{\fg''}
\cong \fh_{r|s}\oplus\fh_{r'|s'},\quad\mbox{and}\quad
\fu_\fq\cong\fu_{\fg'}\oplus\fu_{\fg''}.\]
By \eqref{eq_qcongglrsglrpsp}, all the irreducible polynomial representations of $\fq$ are of the form $L^D_{r|s}\otimes L^E_{r'|s'}$ for some $D\in\Lambda_{r|s}$ and $E\in\Lambda_{r'|s'}$.

In this section, we consider the following branching problem:
  
 \medskip\noindent
 {\bf
 Branching Problem:} {\em  For $F\in\Lambda^+_{p|q}$,     determine a decomposition of the irreducible $\fg$-module $L^F_{p|q}$ into irreducible $\fq$-modules.  
 }

A solution of this branching problem was alluded to in \cite[Remark 9.4]{CZ}. 
Here we want to build a branching algebra to solve the problem.    
  
  \subsection{The branching algebra $\calq(\fg,\fq)$}
 Let $n$ be a positive integer.
 We again consider the $\gl_n\oplus \fg$ module  $\calr=\bigs(\C^n\otimes\C^{p|q})$ 
  and 
  let 
  \[\calq(\fg,\fq)=\calr^{\fu_n\oplus \fu_{\fq}}\] be the subalgebra of $\calr$ consisting of vectors annihilated by the operators from  $\fu_n\oplus \fu_\fq$ (recall that  $\fu_n$
 is the subalgebra of $\gl_n$ consisting of all $n\times n$ strictly upper triangular matrices). 
 Then $\calq(\fg,\fq)$ is a module for $\fh_n\oplus\fh_{\fq}$.
  Since  $\fh_n\oplus\fh_{\fq}$ is isomorphic to $ \fh_n\oplus\fh_{r|s}\oplus\fh_{r'|s'}$, 
  we shall instead regard  $\calq(\fg,\fq)$ as a module for $ \fh_n\oplus\fh_{r|s}\oplus\fh_{r'|s'}$.
  
 Let  
 \[
  \Omega(\fg,\fq)=\Lambda^+_{n,p|q}\times\Lambda^+_{n,r|s}\times \Lambda^+_{n,r'|s'}.
 \]
For each  $(F,D,E)\in \Omega(\fg,\fq)$, 
let $\psi^{(F,D,E)}:=(\psi^F_n, \psi^{D^\sharp}_{r|s}, \psi^{E^\sharp}_{r'|s'})$, 
which is the linear functional of $ \fh_n\oplus\fh_{r|s}\oplus\fkh_{r'|s'}$ defined by equation \eqref{notation_psi}.
Define 
\[
\calq(\fg,\fq)_{(F,D,E)}=\left\{v\in\calq(\fg,\fq): t.v=
   \psi^{(F,D,E)}(t)v,  \   \forall t= (t_1,t_2,t_3)\in\fh_n\oplus\fh_{r|s}\oplus\fkh_{r'|s'}\right\}.
\]

We have the following result.

\begin{proposition}\label{prop_qgg0decom}
\begin{enumerate}
\item[(a)]
The algebra $\calq(\fg,\fq)$ has a direct sum decomposition given by 
\begin{equation}\label{eq_qgqisodecom}
\calq(\fg,\fq)=\bigoplus_{(F,D,E)\in\Omega(\fg,\fq)}\calq(\fg,\fq)_{(F,D,E)}.
\end{equation}

\item[(b)]  For all $(F,D,E)\in\Omega(\fg,\fq)$,     
\begin{equation}\dim \calq(\fg,\fq)_{(F,D,E)}
 =\dim \Hom_{\gl_n}(\rho^F_n,\rho^D_n\otimes\rho^{E}_n)  .\label{eq_dimqgq}
 \end{equation}

 \item[(c)]  For all $(F,D,E)\in\Omega(\fg,\fq)$,     
 \begin{equation}
\dim \calq(\fg,\fg_\fq)_{(F,D,E)}
=\dim\Hom_{\fq}(L^D_{r|s}\otimes\rho^E_{r'|s'}, L^F_{p|q}).\label{eq_dimqgq1}
\end{equation}

\end{enumerate}
\end{proposition}

\begin{proof}
The following isomorphism of algebras is clear, 
\begin{align*}
  \calr   &\cong\bigs(\C^n\otimes\C^{r|s})\otimes\bigs(\C^n\otimes\C^{r'|s'}),
  \end{align*}
which can also be easily proved using arguments similar to those in the proof of Proposition \ref{prop_qgm2}. 
 Using  the $(\gl_n,\gl_{r|s})$-duality  and  $(\gl_n,\gl_{r'|s'})$-duality, we obtain
  \begin{align*}
  \calr&\cong\left(\bigoplus_{D\in\Lambda^+_{n,r|s}} \rho^D_n\otimes L^D_{r|s}\right)\otimes
  \left(\bigoplus_{E\in \Lambda_{n,r'|s'}} \rho^{E}_n\otimes L^E_{r'|s'}\right)\\
 & \cong\bigoplus_{(D,E)\in\Lambda^+_{n,r|s}\times \Lambda^+_{n,r'|s'}} (\rho^D_n\otimes\rho^{E}_n)\otimes L^D_{r|s} \otimes L^{E}_{r'|s'}.
  \end{align*}
  By extracting the vectors in $\calr$ which are annihilated by the operators from $\fu_n\oplus \fu_\fq$, we obtain 
  \begin{equation}\label{eq_qgqdecom1}
  \calq(\fg,\fq)
   \cong\bigoplus_{(D,E)\in\Lambda^+_{n,r|s}\times \Lambda^+_{n,r'|s'}} (\rho^D_n\otimes\rho^{E}_n)^{\fu_n}\otimes (L^D_{r|s})^{\fu_{r|s}} \otimes 
  (L^{E}_{r'|s'})^{\fu_{r'|s'}}.
  \end{equation}
  Let $(D,E)\in\Lambda^+_{n,r|s}\times \Lambda^+_{n,r'|s'}$. Then the space  $(\rho^D_n\otimes\rho^{E}_n)^{\fu_n}$ of vectors in  $\rho^D_n\otimes\rho^{E}_n$ annihilated by the operators from $\fu_n$ is a module for $\fh_n$, and  can be decomposed as  
  \begin{equation}\label{eq_rhodnrhoendecom}(\rho^D_n\otimes\rho^{E}_n)^{\fu_n}=\bigoplus_F(\rho^D_n\otimes\rho^{E}_n)^{\fu_n}_F\end{equation}
  where the direct sum is taken over a set of Young diagrams $F$  and  
 \begin{equation}\label{eq_rhodnrhoef}
  (\rho^D_n\otimes\rho^{E}_n)^{\fu_n}_F=\left\{v\in (\rho^D_n\otimes\rho^{E})^{\fu_n}:\ 
   t.v=\psi^F_n(t)v\ \forall t\in\fu_n\right\},
\end{equation}
i.e., $(\rho^D_n\otimes\rho^{E}_n)^{\fu_n}_F$ is the space of $\gl_n$ highest weight vectors of weight $\psi^F_n$ in $\rho^D_n\otimes\rho^{E}_n$. Therefore, 
\begin{equation}\label{eq_rhodef2101}
\dim (\rho^D_n\otimes\rho^{E}_n)^{\fu_n}_F=\dim \Hom_{\gl_n}(\rho^F_n,\rho^D_n\otimes\rho^{E}_n) ,
\end{equation}    the multiplicity of $\rho^F_n$ in the tensor product $\rho^D_n\otimes\rho^{E}_n$.  
By Theorem \ref{thm:duality-super}, the irreducible representations of $\gl_n$ which occur in $\calr$ are labeled by the Young digrams in $\Lambda^+_{n,p|q}$. Thus, we may assume that the Young diagrams $F$ which appear in the direct sum \eqref{eq_rhodnrhoendecom} also belong to $\Lambda^+_{n,p|q}$. 

By combining equations \eqref{eq_qgqdecom1} and \eqref{eq_rhodnrhoendecom}, we obtain
  \begin{equation}\label{eq_qgqdecom2} 
  \calq(\fg,\fq) 
    \cong\bigoplus_{(F,D,E)\in\Omega(\fg,\fq)}
     (\rho^D_n\otimes\rho^{E}_n)^{\fu_n}_F\otimes(L^D_{r|s})^{\fu_{r|s}} \otimes(L^{E}_{r'|s'})^{\fu_{r'|s'}}.
 \end{equation}
 Now observe that  for all $(F,D,E)\in\Omega(\fg,\fq)$,    
 \begin{equation}\label{eq_qgqgln}
 \calq(\fg,\fq)_{(F,D,E)}\cong 
  (\rho^D_n\otimes\rho^{E}_n)^{\fu_n}_F\otimes(L^D_{r|s})^{\fu_{r|s}} \otimes(L^{E}_{r'|s'})^{\fu_{r'|s'}}
 \end{equation}
 as a module for $\fh_n\oplus\fh_p\oplus\fh_q$.  
 This immediately leads to \eqref{eq_qgqisodecom}. 
 Moreover, by equations \eqref{eq_qgqgln} and \eqref{eq_rhodef2101}, 
  \begin{align*}
  \dim \calq(\fg,\fq)_{(F,D,E)}
 &=\dim \left((\rho^D_n\otimes\rho^{E}_n)^{\fu_n}_F\otimes(L^D_{r|s})^{\fu_{r|s}} \otimes(L^{E}_{r'|s'})^{\fu_{r'|s'}}\right)\\
  &=\dim (\rho^D_n\otimes\rho^{E}_n)^{\fu_n}_F\ \dim( L^D_{r|s})^{\fu_{r|s}} \ \dim(L^{E}_{r'|s'})^{\fu_{r'|s'}}\\
&=\dim (\rho^D_n\otimes\rho^{E}_n)^{\fu_n}_F\\
&=\dim \Hom_{\gl_n}(\rho^F_n,\rho^D_n\otimes\rho^{E}_n) ,
  \end{align*}
 since $\dim( L^D_{r|s})^{\fu_{r|s}} =\dim(L^{E}_{r'|s'})^{\fu_{r'|s'}}=1$. This proves (b).
 

 Next, we use equation \eqref{eq_rdecom} to obtain another decomposition of $\calq(\fg,\fq)$. We have
  \begin{equation}\label{eq_qgqdecom2-1} \calq(\fg,\fq)\cong\bigoplus_{F\in \Lambda^+_{n,p|q}}(\rho^F_n)^{\fu_n}\otimes (L^F_{p|q})^{\fu_{\fq}}
 \end{equation}
 where
$ (L^F_{p|q})^{\fu_{\fq}}$  for each $F\in \Lambda^+_{n,p|q}$  is the space of vectors in  
$L^F_{p|q}$ annihilated by the operators from $\fu_\fq$,  which is a module for $\fh_{\fq} 
\cong \fh_{r|s}\oplus\fh_{r'|s'}$. So $ (L^F_{p|q})^{\fu_{\fq}}$ can be written as a direct sum
 \begin{equation}\label{eq_lfpqug0} 
  (L^F_{p|q})^{\fu_{\fq}}=\bigoplus_{(D,E)\in\Lambda^+_{n,r|s}\times \Lambda^+_{n,r'|s'}}  (L^F_{p|q})^{\fu_{\fq}}_{(D,E)} \end{equation}
where  
 \[
 (L^F_{p|q})^{\fu_{\fq}}_{(D,E)}
=\left\{v\in(L^F_{p|q})^{\fu_{\fq}}:\  t.v
=(\psi^{D^\sharp}_{r|s}, \psi^{E^\sharp}_{r'|s'})(t)v, \ \forall t=(t_1,t_2)\in\fh_{r|s}\oplus\fh_{r'|s'}\right\},
\]
which is  the space of all $\fq$ highest weight vectors in $L^F_{p|q}$ of weight $(\psi^{D^\sharp}_{r|s}, \psi^{E^\sharp}_{r'|s'})$. Hence
 \begin{equation}\label{eq_dimlfpquq} \dim (L^F_{p|q})^{\fu_{\fq}}_{(D,E)}=\dim\Hom_{\fq}(L^D_{r|s}\otimes L^E_{r'|s'},L^F_{p|q}),\end{equation}
the multiplicity of the irreducible $\fq$ representation $L^D_{r|s}\otimes\rho^E_{r'|s'}$ in the irreducible $\fg$ representation $L^F_{p|q}$.
By combining equations \eqref{eq_qgqdecom2-1}  and \eqref{eq_lfpqug0}, we obtain
 \begin{align}
\calq(\fg,\fq)&\cong\bigoplus_{F\in \Lambda^+_{n,p|q}}(\rho^F_n)^{\fu_n}\otimes \left(\bigoplus_{(D,E)\in\Lambda^+_{\min(n,p)}\times R_{n,q}}(L^F_{p|q})^{\fu_{\fq}}_{(D,E)}\right)\nonumber
\\
&\cong\bigoplus_{(F,D,E)\in\Omega(\fg,\fq)}(\rho^F_n)^{\fu_n}\otimes (L^F_{p|q})^{\fu_{\fq}}_{(D,E)
}.\label{eq_qdecom2}\end{align}
Now for each $(F,D,E)\in\Omega(\fg,\fq)$,     we have 
\[\calq(\fg,\fq)_{(F,D,E)}\cong (\rho^F_n)^{\fu_n}\otimes (L^F_{p|q})^{\fu_{\fq}}_{(D,E)
}\]
as a module for $\fh_n\oplus\fh_{r|s}\oplus\fh_{r'|s'}$. 
Since $\dim (\rho^F_n)^{\fu_n}=1$, by equation \eqref{eq_dimlfpquq}, we obtain 
\[
\dim \calq(\fg,\fq)_{(F,D,E)}
 =\dim  (L^F_{p|q})^{\fu_{\fq}}_{(D,E)}=\dim\Hom_{\fq}(L^D_{r|s}\otimes L^E_{r'|s'},L^F_{p|q}),
\]
which proves \eqref{eq_dimqgq1}. 
This completes the proof of the proposition. 
\end{proof}

 \begin{remark}
 \begin{enumerate}\label{rmk:recipro}
 \item[(i)] Part (c) of 
  Proposition \eqref{prop_qgg0decom} shows that the structure of  the algebra $\calq(\fg,\fq) $ encodes part of the branching rule for $(\fg,\fq)$.  In view of this fact, we call $\calq(\fg,\fq)$ a {\bf\em  branching algebra} for $(\fg,\fq)$.
 \item[(ii)] We see from Part (b) of  Proposition \eqref{prop_qgg0decom} that the structure of $\calq(\fg,\fq)$ also encodes information on the decomposition of the tensor products of certain irreducible polynomial  representations of $\gl_n$, which can be viewed as a branching rule for $(\gl_n\oplus\gl_n,\Delta(\gl_n))$,
 where $\Delta(\gl_n)=\Delta^{(1)}(\gl_n)$ in the notation of Remark \ref{rmk:co-mult}. 
 Therefore,  $\calq(\fg,\fq)$ is also a {\bf\em branching algebra for 
 $(\gl_n\oplus\gl_n,\Delta(\gl_n))$.}
 \item[(iii)] Since the structure of  the algebra $\calq(\fg,\fq) $ encodes two sets of branching rules connected by a reciprocity law, it is a {\bf\em reciprocity algebra} in the sense of \cite{HTW1,HTW2}.
 \end{enumerate}
 \end{remark}
 
 \subsection{A branching rule for $(\fg,\fq)$}
 We now derive a branching rule for $(\fg,\fq)$ using the reciprocity algebra $\calq(\fg,\fq)$, 
 and the Littlewood-Richardson rule which we now recall.
  For  a skew tableau $T$, the
{\bf\em word} of $T$ is the sequence $w(T)$ of positive integers
obtained by reading the entries of $T$ from  top to bottom, and
right to left in each row. 
 A {\bf\em Littlewood-Richardson}  tableau is a skew tableau
$T$ which is semistandard and  satisfies the  Yamanouchi word condition,  
that is,
for each positive integer $j$,  starting from the first entry of
$w(T)$ to any place in $w(T)$, there are at least as many $j$s as
$(j+1)$s.
The {\bf\em  Littlewood-Richardson rule} \cite{Fu, HL} states that  if $D$,$E$ and $F$ are Young diagrams with at most $n$ rows, then the
multiplicity $c^F_{D,E}$ of $\rho^F_n$ in the tensor product $\rho^D_n\otimes
\rho^E_n$  is equal the  number of Littlewood-Richardson tableaux
of shape $F/D$ and content $E$.

\begin{theorem}\label{thm_brgq} {\rm\bf(Branching from $\gl_{p|q}$ to $\gl_{r|s}\oplus\gl_{r'|s'}$)} For any $F\in\Lambda^+_{n,p|q}$, 
  \[L^F_{p|q}=\bigoplus_{(D,E)\in\Lambda^+_{n,r|s}\times \Lambda^+_{n,r'|s'}}c^F_{D,E}\ L^D_{r|s}\otimes L^E_{r'|s'}\]
  as a representation for $\fq$.
  \end{theorem}
  \begin{proof} Let $F\in\Lambda^+_{p|q}$. For any $(D,E)\in\Lambda^+_{n,r|s}\times \Lambda^+_{n,r'|s'}$, we have
  \begin{align*}
  \dim\Hom_{\fq}(L^D_{r|s}\otimes L^E_{r'|s'}, L^F_{p|q})&=
  \dim \calq(\fg,\fq)_{(F,D,E)} 
=\dim \Hom_{\gl_n}(\rho^F_n,\rho^D_n\otimes\rho^{E}_n),
\end{align*}
by (b) and (c) of Proposition \ref{prop_qgg0decom}. The right hand side is equal to $c^F_{D,E}$ by the  Littlewood-Richardson rule. Hence 
\[
\dim\Hom_{\fq}(L^D_{r|s}\otimes L^E_{r'|s'}, L^F_{p|q})= c^F_{D,E},
\]
proving the theorem.   
 \end{proof}
  
 \begin{remark}
The branching rules of the oscillator representations of 
the general linear and orthosymplectic Lie superalgebras
were determined in \cite{CLZ, CZ} by exploiting the connection 
between branching and tensor product of representations of dual pairs of Lie (super)algebras, 
but without the framework of branching algebras. 
Also, the branching rule in Theorem \ref{thm_brgq} was described in \cite[Remark 9.4]{CZ}. 
\end{remark}

   We now consider the special case $(r,s,r',s')=(p,0,0,q)$. In this case, $\gl_{r|s}=\gl_p$,
   $\gl_{r'|s'}=\gl_q$ and
   \[\fq=\fg_{\bar0}\cong\gl_p\oplus\gl_q\]
   is the even subspace of $\fg$.

  \begin{corollary}\label{cor_gg0} {\rm\bf(Branching from $\gl_{p|q}$ to $\gl_p\oplus\gl_q$)} For $F\in\Lambda^+_{n,p|q}$, we have
  \begin{equation}\label{eq_gg0}
  L^F_{p|q}=\bigoplus_{(D,E)\in\Lambda^+_{\min(n,p)}\times R_{n,q}}c^F_{D,E^t}\ \rho^D_p\otimes\rho^{E^t}_q\end{equation}
  as a representation for $\fg_{\bar0}$, where $R_{n,q}$ is defined in \eqref{eq_rnqdef}.
  \end{corollary}
  \begin{proof} Let $(r,s,r',s')=(p,0,0,q)$. Then
  \[\Lambda^+_{n,r|s}=\Lambda^+_{n,p|0}=\Lambda^+_{\min(n,p)},\quad 
  \Lambda^+_{n,r'|s'}=\Lambda^+_{n,0|q}=R_{n,q},\]
  and for $(D,E)\in\Lambda^+_{\min(n,p)}\times R_{n,q}$,
  \[L^D_{r|s}=L^D_{p|0}=\rho^D_p,\quad  L^E_{r'|s'}=L^E_{0|q}=\rho^{E^t}_q.\]
  Equation \eqref{eq_gg0} then follows from Theorem \ref{thm_brgq}.
  \end{proof}

 \begin{remark} 
 Let $F=(\lambda_1,\lambda_2,...)\in\Lambda^+_{n,p|q}$. 
 Then the highest weight of $L^F_{p|q}$ is 
 $F^\sharp=(\lambda_1,...,\lambda_p;\mu_1,...,\mu_q)$, 
  where for each $1\le j\le q$, $\mu_j=\max(\lambda^\prime_j-p,0)$.
  Let $D_0$ and $E_0$ be defined respectively by
  \[D_0=(\lambda_1,...,\lambda_p),\quad E_0=(\mu_1,...,\mu_q).\]
  Then by Corollary \ref{cor_gg0}, the $\fg_{\bar0}$ module $\rho^{D_0}_p\otimes\rho^{E_0}_q$ occurs in $L^F_{p|q}$ with multiplicity $1$. Moreover, the $\gl_p\oplus \gl_q$ highest weight vector in  $\rho^{D_0}_p\otimes\rho^{E_0}_q$ is also the $\gl_{p|q}$ highest weight vector in $L^F_{p|q}$
\end{remark}

\subsection{Alternative formulae of branching multiplicities} 
An interesting application of Theorem \ref{thm_brgq} is to deriving alternative formulae for the 
$(\gl_{p|q}, \gl_{r|s}\oplus \fh_{p-r}\oplus\fh_{q-s})$ and $(\gl_{p|q}, \gl_{r|s})$
branching multiplicities given in Theorem \ref{thm:brg-m}.
We briefly discuss it here. 

Let $F\in \Lambda^+_{n,p|q}$. Then by Theorem \ref{thm_brgq}, we have
  \[L^F_{p|q}=\bigoplus_{(D,E)\in\Lambda^+_{n,r|s}\times \Lambda^+_{n,r'|s'}}c^F_{D,E}\ L^D_{r|s}\otimes L^E_{r'|s'}\]
  which can be written as 
  \[
  L^F_{p|q}=\bigoplus_{D\in\Lambda^+_{n,r|s}} \ L^D_{r|s}\otimes\left(\bigoplus_{E\in \Lambda^+_{n,r'|s'}}c^F_{D,E}L^E_{r'|s'}\right).\]
  By writing $L^E_{r'|s'}=\bigoplus_{(\alpha,\beta)\in\Z^{r'}_+\times\Z^{s'}_+}(L^E_{r'|s'})_{(\alpha,\beta)}$, we also have
   \[
  L^F_{p|q}=\bigoplus_{(D,\alpha,\beta)\in\Lambda^+_{n,r|s}\times\Z^{r'}_+\times\Z^{s'}_+} \ L^D_{r|s}\otimes\left(\bigoplus_{E\in \Lambda^+_{n,r'|s'}}c^F_{D,E}(L^E_{r'|s'})_{(\alpha,\beta)}\right).\]
 For any $(D,\alpha,\beta)\in\Lambda^+_{n,r|s}\times\Z^{r'}_+\times\Z^{s'}_+$, 
 \[
 \baln
  \dim\Hom_{\gl_{r|s}}(L^D_{r|s}\otimes\C_{\psi^\alpha_{r'}}\otimes\C_{\psi^\beta_{s'}},L^F_{p|q})
  &=\dim \left(\bigoplus_{E\in \Lambda^+_{n,r'|s'}}c^F_{D,E}(L^E_{r'|s'})_{(\alpha,\beta)}\right)\\
    &=\sum_Ec^F_{D,E}\dim (L^E_{r'|s'})_{(\alpha,\beta)}.
  \ealn
 \]
 Using Corollary \ref{cor:CZ-basis} (i) to the right hand side, we obtain
  \begin{align*}
  \sum_Ec^F_{D,E} N'(E,\alpha,\beta)
  &=\sum_Ec^F_{D,E} \sum_HK_{H,\alpha}K_{E^t/H^t,\beta}
   =\sum_{E,H}c^F_{D,E} K_{H,\alpha}K_{E^t/H^t,\beta},
  \end{align*}
  where the sum is taken over all Young diagrams $E$ and $H$. 
  Hence
  \[
   \dim\Hom_{\gl_{r|s}}(L^D_{r|s}\otimes\C_{\psi^\alpha_{r'}}\otimes\C_{\psi^\beta_{s'}},L^F_{p|q})
   =\sum_{E,H}c^F_{D,E} K_{H,\alpha}K_{E^t/H^t,\beta}.
  \]
    It also follows  that   
  \begin{align*}
  \dim\Hom_{\gl_{r|s}}(L^D_{r|s},L^F_{p|q})
  &=\sum_{(\alpha,\beta)\in\Z^{r'}_+\times\Z^{s'}_+} \dim\Hom_{\gl_{r|s}}(L^D_{r|s}\otimes\C_{\psi^\alpha_{r'}}\otimes\C_{\psi^\beta_{s'}},L^F_{p|q})\\
  &=\sum_{(\alpha,\beta)\in\Z^{r'}_+\times\Z^{s'}_+}\sum_{E,H}c^F_{D,E} K_{H,\alpha}K_{E^t/H^t,\beta}.
  \end{align*}

 \section{Weight vectors of  $L^F_{p|q}$ associated to tableaux}\label{sec3}

 In  Section \ref{sec_glrs}, we obtain a branching rule from $\gl_{p|q}$ to the subalgebra $\fm\cong\gl_{r|s}\oplus\fh_{r'}\oplus \fh_{s'}$.  In this section, we fix $F\in\Lambda^+_{n,p|q}$ and construct explicitly a particular set of $\gl_{p|q}$ weight vectors in $L^F_{p|q}$. We show that when $s=0$ or $s=1$, this set of weight vectors  forms a basis for 
 the space of all $\gl_{r|s}$ highest weight vectors in $L^F_{p|q}$.  In particular, when $r=1$ and $s=0$, this set is the basis of weight vectors for $L^F_{p|q}$ constructed in \cite[Theorem 3.3]{CZ}.


\subsection{Ordered monomials and leading monomials}\label{sec:mo}
 
In this subsection, we define a basis $\calm $ for the algebra $ \calr=\bigs(\C^n\otimes\C^{p|q})$ and use  this basis 
to define the leading monomial of a non-zero element in $\calr$. The notion of leading monomial 
provides an effective way to prove linearly independence of subsets of $\calr$.

Recall from Section \ref{sect:glpq} that $\{e_1, e_2, \dots, e_{p+q}\}$ is the standard basis   for $\C^{p|q}$, 
where  $\{e_1,...,e_p\}$  is the standard basis for the even subspace $(\C^{p|q})_{\bar0}=\C^p$, 
and $\{e_{p+1},...,e_{p+q)}\}$ is the standard basis for the odd subspace $(\C^{p|q})_{\bar1}=\C^q$. 
Let $\{\ve_1,...,\ve_n\}$ be the standard basis for $\C^n$, which is taken to be even. 
For $1\le i\le n$, $1\le j\le p$ and $1\le k\le q$, denote
\[
e_{ij}=\ve_i\otimes e_j,\quad f_{ik}=\ve_i\otimes e_{p+k}, 
\]
and let $\calb_{n,p|q}=\calb_{n,p}\cup \calb_{n,q}$ with 
\[
\calb_{n,p}=\{e_{ij}: 1\le i\le n, 1\le j\le p\},\quad\calb_{n,q}=\{f_{ik}: 1\le i\le n, 1\le k\le q\}.
\]
Then $\calr=\BS(\C^n\ot\C^{p|q})$ is the associative algebra generated by $\calb_{n,p| q}$ with the following relations:
For all $i, i., j, j', k, k'$.
\beq\label{eq:R-ef}
e_{i j} e_{i' j'}= e_{i' j'} e_{i j}, \quad  f_{i k} f_{i' k'} = - f_{i' k'} f_{i k}, \quad e_{i j} f_{i' k'} = f_{i' k'} e_{i j}. 
\eeq
It is a superalgebra with the $\Z_2$-grading defined as follows. Retaining the notation $[v]$ for the $\Z_2$-degree of any homogeneous element $v\in\calr$,  we have $[e_{i j}]=0$ and $[f_{ik}]=1$ for all $i, j, k$. 
The algebra $\calr$ is also $\Z_+$-graded with all elements of $\calb_{n,p|q}$ having degree $1$. 

\begin{remark}
It is evident that the subalgebra of $\calr$ generated by $\calb_{n, p}$ 
subject to the relevant relations in \eqref{eq:R-ef}  is isomorphic to $S(\C^n\ot \C^p)$, 
and the subalgebra generated by $\calb_{n, q}$ is isomorphic to $\Lam(\C^n\otimes\C^q)$. 
The multiplication induces an algebra isomorphism from their tensor product to $\calr$, 
recovering 
$\BS(\C^n\ot\C^{p|q})\cong S(\C^n\otimes\C^p)\ot\Lam(\C^n\otimes\C^q)$. 
\end{remark}

Let us arrange the elements of $\calb_{n,p|q}$ in a rectangular array  as shown below:
\begin{equation}\label{eq_basism}
\begin{array}{cccc|cccc}
e_{11} & e_{12} & \cdots & e_{1p} & f_{11} & f_{12} & \cdots & f_{1q} \\
e_{21} & e_{22} & \cdots & e_{2p} & f_{21} & f_{22} & \cdots & f_{2q} \\
 \vdots & \vdots & & \vdots & \vdots & \vdots &  & \vdots \\
e_{n1} & e_{n2} & \cdots & e_{np} & f_{n1} & f_{n2} & \cdots & f_{nq}
\end{array}.
\end{equation}

As we will consider $\calr$ as a module for the Lie super subalgebra $\fm\cong\gl_{r|s}\oplus\fh_{r'}\oplus \fh_{s'}$ of $\gl_{p|q}$,  it is more convenient to introduce a different set of notation 
for some elements of $\calb_{n,p|q}$. For $1\le i\le n$, $1\le k\le r'$ and $1\le \ell\le s'$, let
\[e^\prime_{ik}=e_{i(r+k)},\quad f^\prime_{i\ell}=f_{i(s+\ell)}.\]
In this new notation, the array \eqref{eq_basism} becomes
\begin{equation}\label{eq_newnotationbasis}
\begin{array}{ccc|ccc|ccc|ccc}
e_{11} & \cdots& e_{1r} &e^\prime_{11}& \cdots & e^\prime_{1r'} & f_{11}&\cdots & f_{1s} &f^\prime_{11}& \cdots &f^\prime_{1s'} \\
e_{21}& \cdots & e_{2r}& e^\prime_{11}& \cdots & e^\prime_{2r'} & f_{21}&\cdots & f_{2s} &f^\prime_{21}& \cdots & f^\prime_{2s'} \\
 \vdots && \vdots &\vdots& & \vdots & \vdots& & \vdots & &\vdots & \vdots \\
e_{n1}& \cdots & e_{nr}& e^\prime_{11}& \cdots & e^\prime_{nr'} & f_{n1}&\cdots & f_{ns} &f^\prime_{n1}& \cdots & f^\prime_{ns'}
\end{array}.
\end{equation}
We now define an ordering on $\calb_{n,p|q}$ (in the notation \eqref{eq_newnotationbasis}) as follows: 
\begin{enumerate}
\item[(O1)] $e_{ab} > e_{cd}$  
if and only if $d >b$, or $d=b$ and $c>a$.
\item[(O2)] $e^\prime_{ab} > e^\prime_{cd}$  
if and only if $d >b$, or $d=b$ and $c>a$.
\item[(O3)] $f_{ab} > f_{cd}$  
if and only if $d >b$, or $d=b$ and $c>a$.
\item[(O4)] $f^\prime_{ab} > f^\prime_{cd}$  
if and only if $d >b$, or $d=b$ and $c>a$.
 \item[(O5)] $e_{ab} > f_{cd}>e^\prime_{ij}>f^\prime_{k\ell}$ for all indices $a,b,c,d,i,j,k,\ell$.
\end{enumerate}

\medskip\begin{definition}
\begin{enumerate}
\item[(i)] We call a non-zero product of the form 
\[
m=h_1   h_2   \cdots   h_s
\]
where $h_i\in\calb_{n,p|q}$ for $1\le i\le s$ 
a {\bf\em monomial} in $ \calr$. 

\item[(ii)] An {\bf\em ordered monomial} in $ \calr$ is a monomial 
$m= h_1   h_2   \cdots   h_s$ such that 
\[h_1 > h_2 > \cdots > h_s. \]  

\item[(iii)] If $m$ is a monomial in $ \calr$, then by reordering its factors if necessary, we obtain 
a unique ordered monomial $[m]$. 
We call $[m]$ {\bf\em the ordered monomial associated with $m$. } Note that $[m]=\pm m$.
\end{enumerate}
\end{definition}

Let $\calm$ denote the set of all ordered monomials   in $\calr$. 
Then $\calm$ is a basis for $ \calr$.
We now extend the ordering (O1-O5) to $\calm$ by the {\em graded lexicographic order}. 
Specifically:
\begin{enumerate}
\item[(O6)] Given ordered monomials $m=h_1   h_2   \cdots   h_s$ and
$m^\prime=h^\prime_1  h^\prime_2   \cdots   h^\prime_{s^\prime}$, we have 
$m>m^\prime$ if and only if either $s>s^\prime$ or $s=s^\prime$ and there exists $1\le u\le s$ such that 
\[h_i=h^\prime_i\ \ \text{for $1\le i\le u-1$, and}\ \ h_u>h^\prime _u.\]
\end{enumerate}

\begin{definition}
 Let $f$ be a non-zero element of $ \calr$. Then $f$ can be written as a linear combination
 \begin{equation}\label{eq:lcf}
f=\sum_{i=1}^rc_im_i
\end{equation}
 where $0\ne c_i\in\C$, $m_i\in \calm$ for $1\le i\le r$ and $m_1>m_2>\cdots >m_r$. Then the {\bf\em leading ordered monomial} of $f$ is defined as
\[\LM(f)=m_1.\]
That is, $\LM(f)$ is the largest ordered monomial which appears in the linear combination \eqref{eq:lcf}.
 \end{definition}

Since $\calm$ is a basis for $\calr$,  it is evident that if $S$ is a subset of $\calr$ and all the elements of $S$ have distinct leading monomials, then $S$ is linearly independent.

\newcommand{\fraka}{\mathfrak{A}}

\subsection{Determinant over $\calr$}
Let $\fraka$ be a complex algebra (not assumed to be commutative), and let 
$\rM_k(\fraka)$ be the space of all $k\times k$ matrices over $\fraka$.
For  $A=(a_{ij})\in\rM_k(\fraka)$, define the {\bf\em determinant of $A$} by
\begin{equation}\label{eqn:defdet}
\det A=\sum_{\sigma\in \Sym_k} \sgn(\sigma) \,
          a_{\sigma(1) 1}  a_{\sigma(2) 2}  \cdots   a_{\sigma(k) k},
          \end{equation}
where $\Sym_k$ is the symmetric group on $\{1,2,...,k\}$. We also write $\det A$ as
\[\det A =\left|\begin{array}{cccc}
a_{11}&a_{12}&\cdots& a_{1k}\\
a_{21}&a_{22}&\cdots& a_{2k}\\
\vdots&\vdots&&\vdots\\
a_{k1}&a_{k2}&\cdots& a_{kk}
\end{array}\right|.
\]
The following are some standard properties of the determinant:
\begin{enumerate}
\item[(D1)]
The function $A \mapsto \det A$ is multi-linear in the rows and in the columns of the matrix. 
\item[(D2)]
If the matrix $B$ is obtained by swapping two rows of $A$, then $\det B=-\det A$. Consequently, if $A$ has two identical rows, then $\det A=0$.
\item[(D3)]
If the matrix $C$ is obtained by adding a multiple of a row of $A$ to another row of $A$, then $\det C=\det A$.  
\end{enumerate}
If we replace row by column in (D2), then the conclusion may not hold in general.
For example, 
if $n\ge 2$ and $\fraka=\calr$, then 
\begin{equation}\label{eq:example}\left|\begin{array}{cc}
 f_{11}&f_{11}\\
 f_{21}&f_{21} 
\end{array}\right|=2 f_{11}f_{21}\ne 0\quad\mbox{in $\calr$.}
\end{equation}

\subsection{Ordered pairs of tableaux} 
In this subsection, we define a set $\calt(F,D,\alpha,\beta)$ of ordered pairs of tableaux whose cardinality is equal to the non-negative integer $N(F,D,\alpha,\beta)$ defined in equation  \eqref{eq_nfdabdefn}.

\begin{definition} Let $E$ be a Young diagram and $T$  a tableau of shape $E$.
\begin{enumerate}
\item[(a)] We let $T^t$ denote the tableau of shape $E^t$ obtained by flipping the boxes and entries of the tableau $T$.

\item[(b)] If $D$ is a Young diagram such that $D$ sits in  $E$, then $T|_D$ shall denote the tableau of shape $D$ obtained by removing all the boxes of  $T$ not in $D$,  and $T/ D$ shall denote the tableau of skew shape $E/D$ obtained by removing all the boxes of $T$ belonging to $D$.
\end{enumerate}
\end{definition}

  \begin{definition} Let $D\in\Lambda^+_{n,r|s}$. We define a tableau
$H_D$ of shape $D$ as follows:
\begin{enumerate}
\item[(i)] If $1\le i\le\min(\ell(D),r)$, the boxes of $D$ in the $i$-th row are all filled with the number $i$. 
\item[(ii)] For $D$ with $\ell(D)>r$ and consisting of $k$ columns, if $r< i\le\ell(D)$ and $1\le j\le k$, then the box of $D$ in the $i$-th row and $j$-th column is filled with the number $j$. 
\end{enumerate}
Denote by $H_D^{\uparrow}$  the tableau which is formed by the first $\min(\ell(D),r)$ rows of $H_D$, and by $H_D^{\downarrow}$ the skew tableau obtained   by removing all boxes  belonging to  $H^\uparrow_D$ from $H_D$.
\end{definition}

Note that in the special case when $\ell(D)\le r$, we have $H_D=H^\uparrow_D$ and $H_D^{\downarrow}$ is the empty tableau.


\begin{example} Let $n=7$, $r=s=2$ and $D=(3,3,2,2,1)$. Then $D\in\Lambda^+_{7,2|2}$, 
\[ 
  H_D=\tableau[s]
{\tf 1 &\tf 1 &\tf 1 \\
\tf 2  &\tf 2  &\tf  2\\
1 & 2\\
  1&2\\
  1},\qquad   H_D^{\uparrow}=\tableau[s]
{\tf 1 &\tf 1 &\tf 1 \\
\tf 2  &\tf 2  &\tf  2 },\qquad H_D^{\downarrow}=\tableau[s]
{  & & \\
 &  &\\1 & 2\\
  1&2\\
  1}.\]
\end{example}

We now fix $(F,D,\alpha,\beta)\in\Omega(\fg,\fm)$ with $N_{(F,D,\alpha,\beta)}\ne 0$.
For each Young diagram $E$, let
\[\SST(E/D,\alpha)=\{T:\ \mbox{$T$ is a semistandard tableau of skew shape $E/D$ and content $\alpha$}\},\] 
\[\SST(F^t/E^t,\beta)=\{T:\ \mbox{$T$ is a semistandard tableau of skew shape $F^t/E^t$ and content $\beta$}\}.\]
Then by the definition of Kostka numbers, we have 
\beq\label{eq:card-k-t}
|\SST(E/D,\alpha)|=K_{E/D,\alpha}, \quad
|\SST(F^t/E^t,\beta)|= K_{F^t/E^t,\beta}. 
\eeq
Let us introduce the following sets:
\begin{align}
\calt(F,D,\alpha,\beta)_E&=\SST(E/D,\alpha)\times\SST(F^t/E^t,\beta), \nonumber\\
\calt(F,D,\alpha,\beta)&=\bigcup_E\calt(F,D,\alpha,\beta)_E, \label{eq_tfdabdefn}
\end{align}
where the union in \eqref{eq_tfdabdefn} is taken over all Young diagrams $E$. Note that if $\calt(F,D,\alpha,\beta)_E$ $\ne\emptyset$, then $D$ sits insides $E$ and $E^t$ sits insides $F^t$. 
Since $D$ and $F$ are fixed, there are only finitely many Young diagrams $E$ for which $\calt(F,D,\alpha,\beta)_E$ are non-empty. Hence, \eqref{eq_tfdabdefn} is actually a finite disjoint  union of nonempty finite sets.

Now  each element $(T_1,T_2)$ of $\calt(F,D,\alpha,\beta)$ gives rise to  a tableau $T_1\ast T_2$ of shape $F$   defined as follows: 

\begin{definition}\label{defn_t1st2}
For $(T_1,T_2)\in \calt(F,D,\alpha,\beta)_E$, $T_1\ast T_2$ shall denote the tableau of shape $F$ such that:
\begin{enumerate}
\item[(i)] $(T_1\ast T_2)|_D=H_D$.
\item[(ii)] $(T_1\ast T_2)|_E/(T_1\ast T_2)|_D=T_1$.
\item[(iii)] $(T_1\ast T_2)^t/E^t=T_2$.
\end{enumerate}
\end{definition}


\begin{example}\label{ex_deltott} Let $n=7$, $r=s=2$, $p=q=4$, $F=(5,4,3,3,3,3,2)$,
$D=(3,3,2,2,1)$, $\alpha=(2,3)$, $\beta=(3,4)$, $E=(3,3,3,3,2,1,1)$
\[T_1=\tableau[s]
{ &  &  \\
  &  & \\
& & \circled{1} \\
  &&\circled{2}\\
&\circled{2}\\ 
\circled{1}\\
\circled{2} }\qquad\mbox{and}\qquad
  T_2= \tableau[s]
{ &  & &&&& \\
  &  &&&&{\bf 1}&{\bf 2} \\
& & &&{\bf 1}&{\bf 2}\\
{\bf 1}&{\bf 2}\\
{\bf 2} }. \]
Then
\[
T_1\ast T_2=\tableau[s]
{\tf 1 &\tf 1  &\tf 1 &{\bf 1}&{\bf 2} \\
 \tf 2 &\tf 2  &\tf 2 &{\bf 2}\\
1&2 &\circled{1} \\
1  &2&\circled{2}\\
1&\circled{2}&{\bf 1}\\ 
\circled{1}&{\bf 1}&{\bf 2}\\
\circled{2}&{\bf 2} }.
\]
\end{example}

\medskip
\begin{lemma} The cardinality of the set $\calt_{(F,D,\alpha,\beta)}$ is $N_{(F,D,\alpha,\beta)}$, i.e.,
\[|\calt_{(F,D,\alpha,\beta)}|=N_{(F,D,\alpha,\beta)}.\]
\end{lemma}
\begin{proof}  
By \eqref{eq:card-k-t}, we have
\begin{align*}
|\calt(F,D,\alpha,\beta)|
&=\sum_E|\calt(F,D,\alpha,\beta)_E|\\
&=\sum_E|\SST(E/D,\alpha)||\SST(F^t/E^t,\beta)|\\
&= \sum_EK_{E/D,\alpha}K_{F^t/E^t,\beta}\\
&=N_{(F,D,\alpha,\beta)},
\end{align*} 
where the last equality follows from the definition of $N_{(F,D,\alpha,\beta)}$ (see \eqref{eq_tfdabdefn}).
\end{proof}

\subsection{Weight vectors associated to ordered pairs of tableaux}\label{sec_opt}

We now fix $(T_1,T_2)\in \calt(F,D,\alpha,\beta)_E$ and let $T=T_1\ast T_2$. Assuming that $T^{(1)},T^{(2)},...$, $T^{(k)}$ are  all the columns of $T$ counted from left to right. Let $1\le j\le k$ and consider the $j$th column $T^{(j)}$ of $T$. We divide $T^{(j)}$ into a maximum of $4$ parts and call it Type $j$ or Type $0$ according to the following:
\scriptsize
\begin{equation}\label{eq_tjdefn}\mbox{\normalsize (Type $j$)}\qquad T^{(j)}=\tableau
{\tf 1\\
\tf 2\\
\tf \vdots\\
 \tf r\\
 j\\
 j\\
 \vdots\\
 j\\
  \circled{$c_1$} \\
\circled{$c_2$}\\
\vdots\\
\circled{$c_u$}\\ 
{\bf d_1}\\
{\bf d_2}\\
{\bf \vdots}\\
{\bf \vdots}\\
{\bf d_v} }\qquad\qquad\mbox{\normalsize or}\qquad\qquad\mbox{\normalsize (Type $0$)}\qquad
T^{(j)}=\tableau
{\tf 1\\
\tf 2\\
\tf \vdots\\
 \tf \ell\\ 
  \circled{$c_1$} \\
\circled{$c_2$}\\
\vdots\\
\circled{$c_u$}\\ 
{\bf d_1}\\
{\bf d_2}\\
{\bf \vdots}\\
{\bf d_v} }
\end{equation}\\  
\normalsize
where the different parts can be described as follows: 
\begin{itemize}
\item
Part 1: the entries $1,2,...,r$ in Type $j$ and the entries $1,2,...,\ell$ where $\ell\le r$ in Type $0$ which come from  $H_D^\uparrow$;
\item
Part 2: a string of  $j$ in Type $j$ which comes from  $H_D^\downarrow$ (this part is missing in Type $0$);
\item
Part 3: the entries $c_1<\cdots <c_u$ come from $T_1$; 
\item
Part 4: the entries $d_1\le\cdots \le d_v$ come from $T_2$.
\end{itemize}
It is possible that any of  part 1, 3, 4 may not appear in $T^{(j)}$.

Assuming that  $T^{(j)}$ has  length $h$, we define $\Delta_{T^{(j)}}$ as  the determinant
\[\Delta_{T^{(j)}}=
\left|
\begin{array}{ccc|ccc|ccc|ccc}
e_{11}&\cdots&e_{1r}&f_{1j}&\cdots&f_{1j}&e^\prime_{1c_1}&\cdots&e^\prime_{1c_u}&f^\prime_{1d_1}&\cdots&f^\prime_{1d_v}\\
e_{21}&\cdots&e_{2r}&f_{2j}&\cdots&f_{2j}&e^\prime_{2c_1}&\cdots&e^\prime_{2c_u}&f^\prime_{2d_1}&\cdots&f^\prime_{2d_v}\\
\vdots &&\vdots &\vdots&&\vdots&\vdots& &\vdots&\vdots&&\vdots\\
e_{h1}&\cdots&e_{hr}&f_{hj}&\cdots&f_{hj}&e^\prime_{hc_1}&\cdots&e^\prime_{hc_u}&f^\prime_{hd_1}&\cdots&f^\prime_{hd_v}
\end{array}
\right|\]
if $T^{(j)}$ is of Type $j$, and 
\[\Delta_{T^{(j)}}=
\left|
\begin{array}{ccc|ccc|ccc}
e_{11}&\cdots&e_{1\ell} &e^\prime_{1c_1}&\cdots&e^\prime_{1c_u}&f^\prime_{1d_1}&\cdots&f^\prime_{1d_v}\\
e_{21}&\cdots&e_{2\ell} &e^\prime_{2c_1}&\cdots&e^\prime_{2c_u}&f^\prime_{2d_1}&\cdots&f^\prime_{2d_v}\\
\vdots &&\vdots &\vdots& &\vdots&\vdots&&\vdots\\
e_{h1}&\cdots&e_{h\ell} &e^\prime_{hc_1}&\cdots&e^\prime_{hc_u}&f^\prime_{hd_1}&\cdots&f^\prime_{hd_v}
\end{array}
\right|\]
if $T^{(j)}$ is of Type $0$.
If some   parts  do not appear in  $T^{(j)}$, then the corresponding columns will be omitted in the determinant.

\begin{lemma}\label{lem:dtj} Let $T^{(j)}$ be as given in equation \eqref{eq_tjdefn}.  
\begin{enumerate}
\item[(i)] The element  $\Delta_{T^{(j)}}$ of $\calr$  is a $\gl_n$ highest weight vector  
and a $\gl_{p|q}$ weight vector.  
\item[(ii)] If $T^{(j)}$ is of Type $0$ or Type $1$, then $\Delta_{T^{(j)}}$ is also a $\fm$ highest weight vector. 
\item[(iii)] The leading ordered monomial of $\Delta_{T^{(j)}}$ is given by 
\begin{equation}\label{eq_lmdt}\LM(\Delta_{T^{(j)}})=
\left\{\begin{array}{ll}
e_{11}\cdots e_{rr}f_{(r+1)j}\cdots f_{\ell j}e^\prime_{(\ell+1)1}
\cdots e^\prime_{(\ell+u)u}f^\prime_{(\ell+u+1)d_1}
\cdots f^\prime_{hd_v}&\mbox{(Type $j$)}\\
e_{11}\cdots e_{\ell\ell} e^\prime_{(\ell+1)1}
\cdots e^\prime_{(\ell+u)u}f^\prime_{(\ell+u+1)d_1}
\cdots f^\prime_{hd_v}&\mbox{(Type $0$).}
\end{array}\right.
\end{equation}
(We assume that there is a total of $\ell$ entries which come from Part 1 and Part 2 of $T^{(j)}$.) 
 \end{enumerate}
\end{lemma}
\begin{proof} Part (i) and part (ii) are clear. For part (iii), 
since  $c_1<c_2<\cdots <c_u$ and $d_1\le d_2\le \cdots d_v$, we obtain by inspection that $\LM(\Delta_{T^{(j)}})$ is equal to the product of the diagonal entries, which is the ordered monomial given in equation \eqref{eq_lmdt}. 
\end{proof}

We now define the element $\Delta_{(T_1,T_2)}$ in $\calr$ by
\begin{equation}\label{eq_defdeltat}
\Delta_{(T_1,T_2)}=\Delta_{T^{(1)}}\Delta_{T^{(2)}}\cdots\Delta_{T^{(k)}}.
\end{equation}
We also define the ordered monomial $m_{(T_1,T_2)}$ as follows: For each box $b$ in $T=T_1\ast T_2$, $R(b)$ shall denote the row  in which the box $b$ lies and $N(b)$ shall denote the number   in the box $b$.  Then $m_{(T_1,T_2)}$ is the ordered monomial associated with 
\begin{equation}\label{eq_defmonomialt}
 \left(\prod_{b\in H_D^\uparrow} e_{R(b)N(b)}\right)
\left(\prod_{b\in H_D^\downarrow} f_{R(b)N(b)}\right)
\left(\prod_{b\in T_1}e^\prime_{R(b)N(b)}\right)\left(\prod_{b\in T_2}f^\prime_{R(b)N(b)}\right).
\end{equation}
Note that $m_{(T_1,T_2)}$ encodes all the information about $(T_1,T_2)$.

 \begin{example} Let $(T_1,T_2)$ be as given in Example \ref{ex_deltott}, and let $T^{(1)},T^{(2)},T^{(3)},T^{(4)},T^{(5)}$ be the columns of the tableau $T=T_1\ast T_2$. Then we have:
\[
T^{(1)}=\tableau[s]{\tf 1\\ \tf 2\\ 1 \\ 1\\ 1\\ \circled{1}\\ \circled{2}},\qquad
\Delta_{T^{(1)}} = \left|\begin{array}{cc|ccc|cc}
e_{11}&e_{12}&f_{11}&f_{11}&f_{11}&e^\prime_{11}&e^\prime_{12}\\
e_{21}&e_{22}&f_{21}&f_{21}&f_{21}&e^\prime_{21}&e^\prime_{22}\\
e_{31}&e_{32}&f_{31}&f_{31}&f_{31}&e^\prime_{31}&e^\prime_{32}\\
e_{41}&e_{42}&f_{41}&f_{41}&f_{41}&e^\prime_{41}&e^\prime_{42}\\
e_{51}&e_{52}&f_{51}&f_{51}&f_{51}&e^\prime_{51}&e^\prime_{52}\\
e_{61}&e_{62}&f_{61}&f_{61}&f_{61}&e^\prime_{61}&e^\prime_{62}\\
e_{71}&e_{72}&f_{71}&f_{71}&f_{71}&e^\prime_{71}&e^\prime_{72}
\end{array}\right|,\]  
\[\LM(\Delta_{T^{(1)}}) =e_{11}e_{22}f_{31}f_{41}f_{51}e^\prime_{61}e^\prime_{72}.\]

\medskip
\[T^{(2)}=\tableau[s]{\tf 1\\ \tf 2\\ 2 \\ 2\\ \circled{2}\\ {\bf 1}\\ {\bf 2}},\qquad
\Delta_{T^{(2)}}= \left|\begin{array}{cc|cc|c|cc}
e_{11}&e_{12}&f_{12}&f_{12}&e^\prime_{12}&f^\prime_{11}&f^\prime_{12}\\
e_{21}&e_{22}&f_{22}&f_{22}&e^\prime_{22}&f^\prime_{21}&f^\prime_{22}\\
e_{31}&e_{32}&f_{32}&f_{32}&e^\prime_{32}&f^\prime_{31}&f^\prime_{32}\\
e_{41}&e_{42}&f_{42}&f_{42}&e^\prime_{42}&f^\prime_{41}&f^\prime_{42}\\
e_{51}&e_{52}&f_{52}&f_{52}&e^\prime_{52}&f^\prime_{51}&f^\prime_{52}\\
e_{61}&e_{62}&f_{62}&f_{62}&e^\prime_{62}&f^\prime_{61}&f^\prime_{62}\\
e_{71}&e_{72}&f_{72}&f_{72}&e^\prime_{72}&f^\prime_{71}&f^\prime_{72}
\end{array}\right|,
\]
\[\LM(\Delta_{T^{(2)}}) =e_{11}e_{22}f_{32}f_{42}e^\prime_{52}f^\prime_{61}f^\prime_{72}.\]

\medskip
\[T^{(3)}=\tableau[s]{\tf 1\\ \tf 2\\ \circled{1}\\ \circled{2}\\ {\bf 1}\\ {\bf 2}},\qquad
\Delta_{T^{(3)}}= \left|\begin{array}{cc|cc|cc}
e_{11}&e_{12}&e^\prime_{11}&e^\prime_{12}&f^\prime_{11}&f^\prime_{12}\\
e_{21}&e_{22}&e^\prime_{21}&e^\prime_{22}&f^\prime_{21}&f^\prime_{22}\\
e_{31}&e_{32}&e^\prime_{31}&e^\prime_{32}&f^\prime_{31}&f^\prime_{32}\\
e_{41}&e_{42}&e^\prime_{41}&e^\prime_{42}&f^\prime_{41}&f^\prime_{42}\\
e_{51}&e_{52}&e^\prime_{51}&e^\prime_{52}&f^\prime_{51}&f^\prime_{52}\\
e_{61}&e_{62}&e^\prime_{61}&e^\prime_{62}&f^\prime_{61}&f^\prime_{62} 
\end{array}\right|,
\]
\[\LM(\Delta_{T^{(3)}}) =e_{11}e_{22}e^\prime_{31}e^\prime_{42}f^\prime_{51}f^\prime_{62}.\]

\medskip
\[T^{(4)}=\tableau[s]{  {\bf 1}\\ {\bf 2}},\qquad
\Delta_{T^{(4)}}= \left|\begin{array}{cc}
 f^\prime_{11}&f^\prime_{12}\\
 f^\prime_{21}&f^\prime_{22}
 \end{array}\right|,\qquad
\LM(\Delta_{T^{(4)}}) = f^\prime_{11}f^\prime_{22}.\]

\medskip
\[T^{(5)}=\tableau[s]{   {\bf 2}},\qquad
\Delta_{T^{(5)}}=  \LM(\Delta_{T^{(5)}}) =f^\prime_{12},
\]

\medskip\[\Delta_{(T_1,T_2)}=\Delta_{T^{(1)}}\Delta_{T^{(2)}}\Delta_{T^{(3)}}\Delta_{T^{(4)}}\Delta_{T^{(5)}}\]
and
\[\LM(\Delta_{(T_1,T_2)})=e^3_{11}e^3_{22}f_{31}f_{41}f_{51}f_{32}f_{42}e^\prime_{31}e^\prime_{61}e^\prime_{42}e^\prime_{52}e^\prime_{72} f^\prime_{11}f^\prime_{51} f^\prime_{61}f^\prime_{12} f^\prime_{22}f^\prime_{62} f^\prime_{72}.\]
\end{example}

 We will see later that for 
 $(T_1,T_2)\in \calt(F,D,\alpha,\beta)$,
 $\Delta_{(T_1,T_2)}$ is a $\gl_{p|q}$ weight vector in $\calr$. To specify the $\gl_{p|q}$  weight of  $\Delta_{(T_1,T_2)}$, we shall adopt the following notation: 
Since $\fh_{p|q}\cong\fh_{r|s}\oplus\fh_{r'}\oplus\fh_{s'}$, any linear functional $\psi_{p|q}$ of $\fh_{p|q}$ can be identified with a linear functional
 \[(\psi_{r|s},\psi_{r'},\psi_{s'})\]
 of $ \fh_{r|s}\oplus\fh_{r'}\oplus\fh_{s'}$ (defined in equation \eqref{notation_psi}) where 
 $\psi_{r|s}$, $\psi_{r'}$ and $\psi_{s'}$ are some linear functionals of 
 $ \fh_{r|s}$, $\fh_{r'}$ and $\fh_{s'}$ respectively. 
   
\begin{lemma}\label{lem:dtmt} Let $(T_1,T_2)\in \calt(F,D,\alpha,\beta)$. Then we have the following:
\begin{enumerate}
\item[(i)] $\Delta_{(T_1,T_2)}$ is a $\gl_n$ highest weight vector of weight $\psi^F_n$
and a $\gl_{p|q}$ weight vector of weight  $(\psi^{D^\sharp}_{r|s},\psi^\alpha_{r'},\psi^\beta_{s'})$.

\item[(ii)] $\LM(\Delta_{(T_1,T_2)})=m_{(T_1,T_2)}$.
 \end{enumerate}
\end{lemma}
\begin{proof} Part (i) is clear. For Part (ii), we have 
\begin{align*}
\LM(\Delta_{(T_1,T_2)})&=\LM(\Delta_{T^{(1)}}\Delta_{T^{(2)}}\cdots\Delta_{T^{(k)}})\\
&=[\LM(\Delta_{T^{(1)}})\LM(\Delta_{T^{(2)}})\cdots\LM(\Delta_{T^{(k)}})]\\
&=m_{(T_1,T_2)}
\end{align*}
by equation \eqref{eq_lmdt}.
\end{proof}
 
 
\begin{corollary}\label{cor_bli}
Let $(F,D,\alpha,\beta)\in\Omega$ and 
\[
\bB(F,D,\alpha,\beta)=\{\Delta_{(T_1,T_2)}:\ {(T_1,T_2)} \in\calt(F,D,\alpha,\beta)\}.
\]
Then $\bB(F,D,\alpha,\beta)$ is a linearly independent set of $\gl_{p|q}$ weight vectors in $L^F_{p|q}$ with weight $(\psi^{D^\sharp}_{r|s},\psi^\alpha_{r'},\psi^\beta_{s'})$.

\end{corollary}
\begin{proof} By Lemma \ref{lem:dtmt}, the elements of $\bB(F,D,\alpha,\beta)$ have distinct leading monomials. Hence, $\bB(F,D,\alpha,\beta)$ is linearly independent.  
\end{proof}


 In the remainder of this section, we shall discuss cases when
 $\bB(F,D,\alpha,\beta)$ yields $\fm$ highest weight vectors. 

\subsection{The $\gl_r$ highest weight vectors in $L^F_{p|q}$}
 In this subsection, we assume that $s=0$.
 In the notation of Section \ref{sec_glrs}, we have
 \[ \fg'\cong \gl_{r|0}=\gl_r, \quad \fm\cong\fg'\oplus \fh_{r'}\oplus\fh_q.\]
 In particular, $\fg'$ is a Lie algebra. 
 The algebra $\calq(\fg,\fm)$ is a module for $\fh_n\oplus\fh_r\oplus\fh_{r'}\oplus\fh_q$ and can be decomposed as 
 \[ \calq(\fg,\fm)=\bigoplus_{(F,D,\alpha,\beta)\in\Omega(\fg,\fm)}\calq(\fg,\fm)_{(F,D,\alpha,\beta)}\]
where 
\[ 
 \Omega(\fg,\fm):=\Lambda^{+}_{n,p|q}\times\Lambda^+_{\min(n,r)}\times \Z^{r'}_{+}\times \Z^{q}_{+}.\]
 Note that the component $\Lambda^+_{n,r|s}$ of $\Omega(\fg,\fm)$ reduces to $\Lambda^+_{\min(n,r)}$ in this case. 
 
 
\begin{theorem}\label{thm_cases0}  
\begin{enumerate}
\item[(i)] If $(F,D,\alpha,\beta)\in\Omega(\fg,\fm)$, then the set
\[
\bB(F,D,\alpha,\beta)=\{\Delta_{(T_1,T_2)}:\ {(T_1,T_2)} \in\calt(F,D,\alpha,\beta)\}
\]
is a basis for $\calq(\fg,\fm)_{(F,D,\alpha,\beta)}$.

\item[(ii)] For $F\in\Lambda^+_{n,p|q}$ and $D\in \Lambda^+_{\min(n,r)}$, the set 
\[\bigcup_{(\alpha,\beta)\in\Z^{r'}_{+}\times\Z^q_{+}}
\bB(F,D,\alpha,\beta)\]
is a basis for the space of $\gl_r$ highest weight vectors of weight $\psi^D_r$ in $L^F_{p|q}$.

\item[(iii)] The set 
\[\bB=\bigcup_{(F,D,\alpha,\beta)\in\Omega}\bB(F,D,\alpha,\beta)\]
is a basis for $\calq(\fg,\fm)$.
\end{enumerate}
\end{theorem}
\begin{proof} Let $(F,D,\alpha,\beta)\in\Omega(\fg,\fm)$, $(T_1,T_2) \in\calt(F,D,\alpha,\beta)$ and $T=T_1\ast T_2$. We recall the definition of $\Delta_{(T_1,T_2)}$ as given in equation \eqref{eq_defdeltat} and assume that the tableau $T$ has $k$ columns. Since $D\in\Lambda^+_{\min(n,r)}$, $\ell(D)\le r$. Consequently, for each $1\le j\le k$, the column tableau $T^{(j)}$ is of Type $0$ and so by Part (ii) of Lemma \ref{eq_tjdefn}, $\Delta_{T^{(j)}}$ is a $\fm$ highest weight vector. Hence,
$\Delta_{(T_1,T_2)}$ is a product of $\fm$ highest weight vectors, and  it is also a $\fm$ highest weight vector.  It follows that $\Delta_{(T_1,T_2)}\in\calq(F,D,\alpha,\beta)$. This proves that $\bB(F,D,\alpha,\beta)\subseteq \calq(F,D,\alpha,\beta)$. 


Since $|\bB(F,D,\alpha,\beta)|=|\calt(F,D,\alpha,\beta)|=\dim \calq(F,D,\alpha,\beta)$ and by Corollary \ref{cor_bli},  $\bB(F,D,\alpha,\beta)$ is linearly independent, it is a basis for 
$\calq(F,D,\alpha,\beta)$. This proves (i).


As explained in the proof of Proposition \ref{prop_qgmdecom}, $\calq(\fg,\fm)_{(F,D,\alpha,\beta)}$ can be identified with the space of $\fm$ highest weight vectors of weight $(\psi^D_r,\psi^\alpha_{r'},\psi^\beta_q)$ in $L^F_{p|q}$. By taking union over all possible $(\alpha,\beta)$, we obtain all the $\gl_r$ highest weight vectors of weight $\psi^D_r$ in $L^F_{p|q}$. This gives (ii).


Finally,
(iii) follows from (i) and equation \eqref{eq_qgmisodecom}.
\end{proof}

\subsection{A basis for $L_{p|q}^F$}
 We now consider the case with $r=s=0$ and will use the notation defined in \S \ref{sec_wtm}. In particular, we recall that in this case,  $\Omega(\fg,\fm)=\Lambda^+_{n,p|q}\times\Z^p_{+}\times\Z^q_{+}$. 
For $(F,\alpha,\beta)\in\Omega(\fg,\fm)$, let  
\[\calt'(F,\alpha,\beta)=\bigcup_E\calt'(F,\alpha,\beta)_E\]
where the union is taken over all  Young diagram $E$, and 
\[\calt'(F,\alpha,\beta)_E=\SST(E,\alpha)\times\SST(F^t/E^t,\beta),\]
that is,   $\calt'(F,\alpha,\beta)_E$    consists of all ordered pairs $(T_1,T_2)$ of tableaux such that $T_1$ is semistandard tableau of shape $E$  and $T_2$ is a semistandard tableau of skew shape $F^t/E^t$.   Then we have
$|\calt'(F,\alpha,\beta)| =N'_{(F,\alpha,\beta)}$ where $N'_{(F,\alpha,\beta)}$ is defined in 
equation \eqref{eq:mult-formula}.
For $(T_1,T_2)\in\calt'(F,\alpha,\beta)_E$, let $T_1\ast T_2$ be the tableau of shape $F$ such that:
\begin{enumerate}
\item[(i)]   $(T_1\ast T_2)|_E=T_1$.
\item[(ii)] $(T_1\ast T_2)^t/E^t=T_2$.
\end{enumerate}

\begin{remark} In \cite{BR},
 $T_1\ast T_2$ is called a {\bf $\mathbf{(p,q)}$ semistandard tableau} of shape $F$. It is also shown there that the dimension of $L^F_{p|q}$ is equal to the number of all such tableaux, which is $\tilde{N}'_F$. See Part (ii) of Corollary \ref{cor:CZ-basis}.
\end{remark}

 We now fix $(T_1,T_2)\in \calt'(F,\alpha,\beta)_E$ and  construct an element $\Delta_{(T_1,T_2)}$ of $\calr$ using a procedure similar to that given in \S \ref{sec_opt}.
 Let $T=T_1\ast T_2$ and assuming that $T^{(1)},T^{(2)},...$, $T^{(k)}$ are  all the columns of $T$ counted from left to right. Let $1\le j\le k$ and consider the $j$th column $T^{(j)}$ of $T$. We divide $T^{(j)}$ into a maximum of $2$ parts as follows: 
\scriptsize
\[  T^{(j)}=\tableau
{ 
  \circled{$c_1$} \\
\circled{$c_2$}\\
\vdots\\
\circled{$c_u$}\\ 
{\bf d_1}\\
{\bf d_2}\\
{\bf \vdots}\\
{\bf \vdots}\\
{\bf d_v} } 
\]  
\normalsize
where   the entries $c_1<\cdots <c_u$ come from $T_1$ and 
 the entries $d_1\le\cdots \le d_v$ come from $T_2$. Let $\Delta_{T^{(j)}}$
 be the determinant  
\[\Delta_{T^{(j)}}=
\left|
\begin{array}{cccccccc}
e_{1c_1}&e_{1c_2}&\cdots&e_{1c_u}&f_{1d_1}&f_{1d_2}&\cdots&f_{1d_v}\\
e_{2c_1}&e_{2c_2}&\cdots&e_{2c_u}&f_{2d_1}&f_{2d_2}&\cdots&f_{2d_v}\\
\vdots&\vdots& \ddots &\vdots&\vdots&\vdots& \ddots &\vdots\\
e_{hc_1}&e_{hc_2}&\cdots&e_{hc_u}&f_{hd_1}&f_{hd_2}&\cdots&f_{hd_v}\\
\end{array}\right|\]
where $h=u+v$.  Then  $\Delta_{(T_1,T_2)}$ is defined as
\[
\Delta_{(T_1,T_2)}=\Delta_{T^{(1)}}\Delta_{T^{(2)}}\cdots\Delta_{T^{(k)}}.
\]

 \begin{corollary} \label{cor:CZ-basis2}
 Let $F\in\Lambda^+_{n,p|q}$.
\begin{enumerate}
\item[(i)] For $(\alpha,\beta)\in\Z^p_{+}\times\Z^q_{+}$,     the set
\[
\bB(F,\alpha,\beta)=\{\Delta_{(T_1,T_2)}:\ {(T_1,T_2)} \in\calt'(F,\alpha,\beta)\}
\]
is a basis for $(L^F_{p|q})_{(\alpha,\beta)}$ (which is identified with $\calq(\fg,\fm)_{(F,\alpha,\beta)}$ by Corollary \ref{cor:CZ-basis} (i)).

\item[(ii)] The set 
\[\bigcup_{(\alpha,\beta)\in\Z^p_{+}\times\Z^q_{+}}
\bB(F,\alpha,\beta)\]
is a basis for   $L^F_{p|q}$ (which is identified with the subspace $(\rho^F_n)^{\fu_n}\otimes L^F_{p|q}$ of $\calq(\fg,\fm)$).

\item[(iii)] The set 
\[\bB=\bigcup_{(F,\alpha,\beta)\in\Omega(\fg,\fm)}\bB(F,\alpha,\beta)\]
is a basis for $\calq(\fg,\fm)$.
\end{enumerate}
 
 \end{corollary}

\begin{remark} 
The basis for $L^F_{p|q}$ given in part (ii) of Corollary \ref{cor:CZ-basis} was obtained in \cite[Theorem 3.3]{CZ} by direct calculations. Here we construct the basis from the iterated Pieri rules (see Section \ref{sec_opt}). This also explains why the semistandard tableaux enter the picture.
\end{remark}

 \subsection{The $\gl_{r|1}$ highest weight vectors in $L^F_{p|q}$}
 In this final subsection, we consider the case $s=1$.
 In this case, we have $\fg'\cong\gl_{r|1}$ and  $\fm\cong\fg'\oplus \fh_{r'}\oplus\fh_{q-1}$.
   The algebra $\calq(\fg,\fm)$ is a module for $\fh_n\oplus\fh_{r|1}\oplus\fh_{r'}\oplus\fh_{q-1}$ 
   and can be decomposed as 
 \[ \calq(\fg,\fm)=\bigoplus_{(F,D,\alpha,\beta)\in\Omega(\fg,\fm)}\calq(\fg,\fm)_{(F,D,\alpha,\beta)}\]
where 
\[ 
 \Omega(\fg,\fm)=\Lambda^{+}_{n,p|q}\times\Lambda^+_{n,r|1}\times \Z^{r'}_{+}\times \Z^{q-1}_{+}.
\]

\begin{theorem}  \label{thm_cases1}  
\begin{enumerate}
\item[(i)] If $(F,D,\alpha,\beta)\in\Omega(\fg,\fm)$, then the set
\[
\bB(F,D,\alpha,\beta)=\{\Delta_{(T_1,T_2)}:\ {(T_1,T_2)} \in\calt(F,D,\alpha,\beta)\}
\]
is a basis for $\calq(\fg,\fm)_{(F,D,\alpha,\beta)}$, which can be identified with the space of all $\fm$ highest weight vectors in $L^F_{p|q}$ of weight $(\psi^{D^\sharp}_{r|1},\psi^\alpha_{r'},\psi^\beta_{q-1})$.

\item[(ii)] For $F\in\Lambda^+_{n,p|q}$ and $D\in \Lambda^+_{n,r|1}$, the set 
\[\bigcup_{(\alpha,\beta)\in\Z^{r'}_{+}\times\Z^q_{+}}
\bB(F,D,\alpha,\beta)\]
is a basis for the space of all $\gl_{r|1}$ highest weight vectors of weight $\psi^{D^\sharp}_{r|1}$ in $L^F_{p|q}$.

\item[(iii)] The set 
\[\bB=\bigcup_{(F,D,\alpha,\beta)\in\Omega(\fg,\fm)}\bB(F,D,\alpha,\beta)\]
is a basis for $\calq(\fg,\fm)$.
\end{enumerate}
\end{theorem}
\begin{proof} 

Let $(F,D,\alpha,\beta)\in\Omega$, $(T_1,T_2) \in\calt(F,D,\alpha,\beta)$ and $T=T_1\ast T_2$. Assume that $T$ has $k$ columns. Then since $D\in\Lambda^+_{r|1}$,  the column tableau $T^{(j)}$ is of Type $1$ for each $1\le j\le k$.  By Part (ii) of Lemma \ref{eq_tjdefn}, each $\Delta_{T^{(j)} }$ is a $\fm$ highest weight vector.
 The remaining arguments for (i) and the  proof for (ii) and (iii) are  similar to that of Theorem \ref{thm_cases0}.  
\end{proof}


\medskip


\end{document}